\pdfoutput=1
\documentclass[pdflatex,preprint,number]{elsarticle}

\usepackage{amssymb}
\usepackage{amsfonts}
\usepackage{hyperref}
\usepackage{epsfig,amsmath}
\usepackage{tikz,ifthen,calc}
\usepackage{graphicx}
\usepackage{caption}
\usepackage{subcaption}
\usepackage{wrapfig}
\usepackage{scalefnt}

\usepackage[margin=1.00 in]{geometry}

\newtheorem{thm}{Theorem}

\newtheorem{propo}[thm]{Proposition}

\newdefinition{rmk}{Remark}
\newdefinition{definition}{Definition}
\newdefinition{example}{Example}
\newproof{pf}{Proof}

\begin{document}

\title{Signed tilings by ribbon $L$ $n$-ominoes, $n$ odd, via Gr\"obner bases}

\author[rvt1]{Viorel Nitica}
\ead{vnitica@wcupa.edu}

\address[rvt1]{Department of Mathematics, West Chester University, PA 19383, USA}

\begin{abstract} We show that a rectangle can be signed tiled by ribbon $L$ $n$-ominoes, $n$ odd, if and only if it has a side divisible by $n$. A consequence of our technique, based on the exhibition of an explicit Gr\"obner basis, is that any $k$-inflated copy of the skewed  $L$ $n$-omino has a signed tiling by skewed $L$ $n$-ominoes. We also discuss regular tilings by ribbon $L$ $n$-ominoes, $n$ odd, for rectangles and more general regions. We show that in this case obstructions appear that are not detected by signed tilings.
\end{abstract}

\begin{keyword}
replicating tile; $L$-shaped polyomino; skewed $L$-shaped polyomino; signed tilings; Gr\"obner basis
\end{keyword}

\date{}
\maketitle

\section{Introduction}

In this article we study tiling problems for regions in a square lattice by certain symmetries of an $L$-shaped polyomino.
Polyominoes were introduced by Solomon W. Golomb in \cite{Golomb1} and the standard reference about this subject is the book \emph{Polyominoes} \cite{Golumb3}. The $L$-shaped polyomino we study is placed in a square lattice and is made out of $n, n\ge 3,$ unit squares, or \emph{cells}. See Figure~\ref{fig:LTetromino1}. In an $a\times b$ rectangle, $a$ is the height and $b$ is the base. We consider translations (only!) of the tiles shown in Figure~\ref{fig:LTetrominoes1}. They are ribbon $L$-shaped $n$-ominoes. A {\em ribbon polyomino}~\cite{Pak} is a simply connected polyomino with no two unit squares lying along a line parallel to the first bisector $y=x$. We denote the set of tiles by $\mathcal{T}_n.$

\begin{figure}[h]
~~~~~~~~~
\begin{subfigure}{.15\textwidth}
\begin{tikzpicture}[scale=.35]
\draw [line width = 1](2.5,1)--(1,1)--(1,2)--(0,2)--(0,0)--(2.5,0);
\draw[dotted, line width=1] (2.6,.5)--(3.9,.5);
\draw [line width = 1] (4,0)--(6,0)--(6,1)--(4,1);
\draw [line width = 1] (1,0)--(1,1);
\draw [line width = 1] (0,1)--(1,1);
\draw [line width = 1] (2,0)--(2,1);
\draw [line width = 1] (5,0)--(5,1);
\end{tikzpicture}
\caption{An $L$ $n$-omino with $n$ cells.}
\label{fig:LTetromino1}
\end{subfigure}
~~~~~~~~~
\begin{subfigure}{.5\textwidth}
\begin{tikzpicture}[scale=.35]
\draw [line width = 1](2.5,1)--(1,1)--(1,2)--(0,2)--(0,0)--(2.5,0);
\draw[dotted, line width=1] (2.6,.5)--(3.9,.5);
\draw [line width = 1] (4,0)--(6,0)--(6,1)--(4,1);
\draw [line width = 1] (1,0)--(1,1);
\draw [line width = 1] (0,1)--(1,1);
\draw [line width = 1] (2,0)--(2,1);
\draw [line width = 1] (5,0)--(5,1);
\node at (3, -1) {$T_3$};
\draw [line width = 1] (10,1)--(8,1)--(8,2)--(10,2);
\draw [dotted, line width=1] (10.1,1.5)--(11.4,1.5);
\draw [line width = 1] (9,1)--(9,2);
\draw [line width = 1] (11.5,1)--(13,1)--(13,0)--(14,0)--(14,2)--(11.5,2);
\draw [line width = 1] (12,1)--(12,2);
\draw [line width = 1] (13,1)--(13,2);
\draw [line width = 1] (13,1)--(14,1);
\node at (11, -1) {$T_4$};

\draw [line width = 1] (-9,2)--(-9,4)--(-8,4)--(-8,2);
\draw [dotted, line width=1] (-8.5,1.9)--(-8.5,.6);
\draw [line width = 1] (-9,.5)--(-9,-2)--(-7,-2)--(-7,-1)--(-8,-1)--(-8,.5);
\draw [line width = 1] (-9,3)--(-8,3);
\draw [line width = 1] (-9,0)--(-8,0);
\draw [line width = 1] (-9,-1)--(-8,-1);
\draw [line width = 1] (-8,-1)--(-8,-2);
\node at (-8, -3) {$T_1$};

\draw [line width = 1] (-3,1.5)--(-3,3)--(-4,3)--(-4,4)--(-2,4)--(-2,1.5);
\draw [line width = 1] (-3,2)--(-2,2);
\draw [line width = 1] (-3,3)--(-2,3);
\draw [line width = 1] (-3,3)--(-3,4);
\draw [dotted, line width=1] (-2.5,1.4)--(-2.5,.1);
\draw [line width = 1] (-3,0)--(-3,-2)--(-2,-2)--(-2,0);
\draw [line width = 1] (-3,-1)--(-2,-1);
\node at (-2.5, -3) {$T_2$};
\end{tikzpicture}
\caption{The set of tiles $\mathcal{T}_n$.}
\label{fig:LTetrominoes1}
\end{subfigure}
\caption{}
\end{figure}

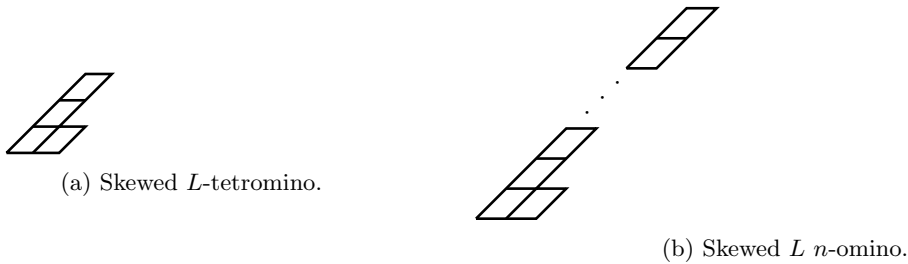
\begin{figure}[h]
~~~~~~~~~
\begin{subfigure}{.30\textwidth}
\begin{tikzpicture}[scale=.35]
\draw [line width = 1](0,0)--(3,3)--(4,3)--(2,1)--(3,1)--(2,0)--(0,0);
\draw [line width = 1](1,1)--(2,1);
\draw [line width = 1](2,2)--(3,2);
\draw [line width = 1](1,0)--(2,1);
\end{tikzpicture}
\caption{Skewed $L$-tetromino.}
\label{fig:skewed-LTetromino1}
\end{subfigure}
~~~~~~~~~
\begin{subfigure}{.5\textwidth}
\begin{tikzpicture}[scale=.4]
\draw [line width = 1](0,0)--(3,3)--(4,3)--(2,1)--(3,1)--(2,0)--(0,0);
\draw [line width = 1](1,1)--(2,1);
\draw [line width = 1](2,2)--(3,2);
\draw [line width = 1](1,0)--(2,1);

\draw [line width = 1] (5,5)--(7,7)--(8,7)--(6,5)--(5,5);
\draw [line width = 1] (6,6)--(7,6);
\node at (3.7, 3.5) {$\cdot$};
\node at (4.2, 4) {$\cdot$};
\node at (4.7, 4.5) {$\cdot$};
\end{tikzpicture}
\caption{Skewed $L$ $n$-omino.}
\label{fig:skewd-LTetrominoesn}
\end{subfigure}
\caption{Skewed polyominoes}
\end{figure}

Tilings by $\mathcal{T}_n, n$ even, are studied in~\cite{CLNS},~\cite{nitica-L-shaped}, with \cite{CLNS} covering the case $n=4$. We recall that a replicating tile is one that can make larger copies of itself. The order of replication is the number of initial tiles that fit in the larger copy. Replicating tiles were introduced by Golomb in~\cite{Golomb-rep}. In~\cite{Nitica} we study replication of higher orders for several replicating tiles introduced in~\cite{Golomb-rep}. In particular, it is suggested there that the skewed $L$-tetromino showed in Figure~\ref{fig:skewed-LTetromino1} is not replicating of order $k^2$ for any odd $k$. The question is equivalent to that of tiling a $k$-inflated copy of the straight $L$-tetromino using only four, out of eight possible, orientations of an $L$-tetromino, namely those that are ribbon. The question is solved in~\cite{CLNS}, where it is shown that the $L$-tetromino is not replicating of any odd order. This is a consequence of a stronger result: a tiling of the first quadrant by $\mathcal{T}_4$ always follows the rectangular pattern, that is, the tiling reduces to a tiling by $4\times 2$ and $2\times 4$ rectangles, each tiled in turn by two tiles from $\mathcal{T}_4$.

The results in~\cite{CLNS} are generalized in~\cite{nitica-L-shaped} to $\mathcal{T}_n, n$ even. The main result shows that any tiling of the first quadrant by $\mathcal{T}_n$ reduces to a tiling by $2\times n$ and $n\times 2$ rectangles, with each rectangle covered by two ribbon $L$-shaped $n$-ominoes. An application is the characterization of all rectangles that can be tiled by $\mathcal{T}_n, n$ even: a rectangle can be tiled if and only if both sides are even and at least one side is divisible by $n$. Another application is the existence of the local move property for an infinite family of sets of tiles: $\mathcal{T}_n, n$ even, has the local move property for the class of rectangular regions with respect to the local moves that interchange a tiling of an $n\times n$ square by $n/2$ vertical  rectangles, with a tiling by $n/2$ horizontal rectangles, each vertical/horizontal rectangle being covered by two ribbon $L$-shaped $n$-ominoes. One shows that neither of these results are valid for any odd $n$. The rectangular pattern of a tiling of the first quadrant persists if one adds an extra $2\times 2$ tile to $\mathcal{T}_n, n$ even. A rectangle can be tiled by the larger set of tiles if and only if it has both sides even. It is also shown in the paper that the main result implies that a skewed $L$-shaped $n$-omino, $n$ even, (see Figure~\ref{fig:skewd-LTetrominoesn}) is not a replicating tile of order $k^2$ for any odd $k$.

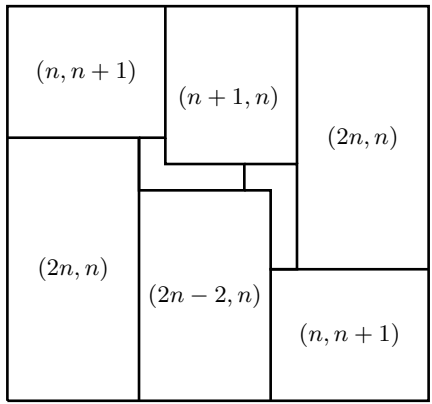
\begin{figure}[h]
\centering
\begin{tikzpicture}[scale=.35]
\draw [line width = 1] (0,0)--(16,0)--(16,15)--(0,15)--(0,0);
\draw [line width = 1] (0,10)--(5,10)--(5,0);
\draw [line width = 1] (5,10)--(6,10)--(6,15);
\draw [line width = 1] (6,10)--(6,9)--(9,9)--(9,8)--(5,8)--(5,10);
\draw [line width = 1] (9,9)--(11,9)--(11,5)--(10,5)--(10,8)--(9,8);
\draw [line width = 1] (10,0)--(10,5)--(16,5);
\draw [line width = 1] (11,9)--(11,15);
\node at (13, 2.5) {\small{$(n,n+1)$}};
\node at (3, 12.5) {\small{$(n,n+1)$}};
\node at (2.5, 5) {\small{$(2n,n)$}};
\node at (13.5, 10) {\small{$(2n,n)$}};
\node at (7.5, 4) {\small{$(2n-2,n)$}};
\node at (8.4, 11.5) {\small{$(n+1,n)$}};
\end{tikzpicture}
\caption{A tiling of a $(3n,3n+1)$ rectangle by $\mathcal{T}_n$.}
\label{fig:odd-rectangle}
\end{figure}

We investigate in this paper tiling properties of $\mathcal{T}_n, n$ odd. Parallel results with~\cite{nitica-L-shaped} are not possible due to the fact, already observed in~\cite{nitica-L-shaped}, that there are rectangles that have tilings by $\mathcal{T}_n, n$ odd, that do not follow the rectangular pattern. See Figure~\ref{fig:odd-rectangle}. Instead of regular tilings one can study signed tilings. These are finite placements of tiles on a plane, with weights +1 or -1 assigned to each of the tiles. We say that they tile a region $R$
if the sum of the weights of the tiles is 1 for every cell inside $R$ and 0 for every cell elsewhere.

A useful tool in the study of signed tilings is a Gr\"obner basis associated to the polynomial ideal generated by the tiling set. If the coordinates of the lower left corner of a cell are $(\alpha,\beta)$, one associates to the cell the monomial $x^{\alpha}y^{\beta}$. This correspondence associates to any bounded tile placed in the square lattice a Laurent polynomial with all coefficients $1$. The polynomial associated to a tile $P$ is denoted by $f_P$. The polynomial associated to a tile translated by an integer vector $(\gamma, \delta)$ is the initial polynomial multiplied by the monomial $x^{\gamma}y^{\delta}$. If the region we want to tile is bounded and if the tile set consists of bounded tiles, then the whole problem can be translated in the first quadrant via a translation by an integer vector, and one can work only with regular polynomials in $\mathbb{Z}[X,Y]$. See Theorem~\ref{thm-groebner-tiling} below.

Our main result is the following:

\begin{thm}\label{thm-main} A rectangle can be signed tiled by $\mathcal{T}_n, n\ge 5$ odd, if and only if it has a side divisible by $n$.
\end{thm}

Theorem~\ref{thm-main} is proved in Section 4 using a Gr\"obner basis for the tiling set computed in Section 3.

For completeness, we briefly discuss regular tilings by $\mathcal{T}_n, n\ge 5$ odd.

Theorem~\ref{thm-main} gives for regular tilings by $\mathcal{T}_n, n\ge 5$ odd, a corollary already known (see~[Lemma 2]\cite{nitica-L-shaped}):

\begin{thm} If $n\ge 5$ odd, a rectangle with neither side divisible by $n$ cannot be tiled by $\mathcal{T}_n$.
\end{thm}

If one of the sides of the rectangle is divisible by $n$, we recall first the following result of Herman Chau, mentioned in~\cite{nitica-L-shaped}, which is based on a deep result of Pak~\cite{Pak}:

\begin{thm} A rectangle with both sides odd cannot be tiled by $\mathcal{T}_n, n\ge 5$ odd.
\end{thm}

If one of the sides of the rectangle is even, one has the following result.

\begin{thm}\label{thm:some445} Let $n\ge 5,$ odd and assume that a rectangle has a side divisible by $n$ and a side of even length.
\begin{enumerate}
\item If one side is divisible by $n$ and the other side is of even length, then the rectangle can be tiled by $\mathcal{T}_n$.
\item If the side divisible by $n$ is of length at least $3n+1$ and even, and the other side is of length at least $3n$ and odd, then the rectangle can be tiled by $\mathcal{T}_n$.
\end{enumerate}
\end{thm}

\begin{pf} 1) The rectangle can be tiled by $2\times n$ or $n\times 2$ rectangles, which can be tiled by two tiles from $\mathcal{T}_n$.

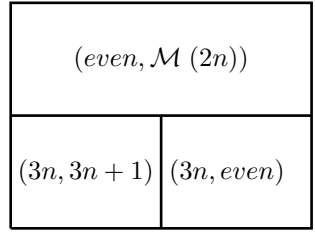
\begin{figure}[h!]
\centering
\begin{tikzpicture}[scale=.25]
\draw [line width = 1] (0,0)--(16,0)--(16,12)--(0,12)--(0,0);
\draw [line width = 1] (8,0)--(8,6)--(0,6);
\draw [line width = 1] (8,6)--(16,6);

\node at (11.5, 3) {$(3n,even)$};
\node at (8, 9) {$(even,{\cal M}\;( 2n))$};
\node at (4, 3) {$(3n,3n+1)$};
\end{tikzpicture}
\caption{A tiling of an (odd, even) rectangle by $\mathcal{T}_n$.}
\label{fig:odd-even-rectangle}
\end{figure}

2) We use the tiling shown in Figure~\ref{fig:odd-even-rectangle}. The $(3n,3n+1)$ rectangle is tiled as in Figure~\ref{fig:odd-rectangle}, and the other two  rectangles can be tiled by $2\times n$ or $n\times 2$ rectangles, which in turn can be tiled by two tiles from $\mathcal{T}_n$.
\end{pf}

A consequence of the technique used in the proof of Theorem~\ref{thm-main} is:

\begin{propo}\label{propo:skewd} If $n\ge 5$ odd, a $k$-inflated copy of the  $L$ $n$-omino has a signed tiling by ribbon $L$ $n$-ominoes.
\end{propo}

Proposition~\ref{propo:skewd} is proved in Section 5.

As any $2n\times 2n$ square can be tiled by $\mathcal{T}_n$, it follows that if $k$ is divisible by $2n$ then the skewed $L$ $n$-omnino is replicating of order $k^2$. Information about other orders of replication can be found using Pak's invariant~\cite{Pak}.

\begin{propo}\label{propo:skew93} Let $n\ge 5$ odd.

1) If $k\ge 1$ is odd and divisible by $n$, then the skewed $L$ $n$-omino is not replicating of order $k^2$.

2) If $k\ge 1$ is even and not divisible by $n$, then the skewed $L$ $n$-omino is not replicating of order $k^2$.
\end{propo}

Proposition~\ref{propo:skew93} is proved in Section 6.

Proposition~\ref{propo:skew93} leaves open the question of replication of the skewed $L$ $n$-omino if $k$ is odd and not divisible by $n$. Some cases can be solved using Pak's higher invariants $f_2, \dots f_m$~\cite{Pak}, which are all zero for tiles in $\mathcal{T}_n$. For example, if $n=5,$ a $3$-inflated copy of the $L$ pentomino has $f_2=-1$, showing the impossibility of tiling.

A general result for regular tilings is out of reach due to the fact that for $k$ odd and congruent to 1 modulo $n$, the leftover region that appears (see the proof of Proposition~\ref{propo:skew93}) is just an $L$ $n$-omino, that has all higher invariants $f_2, \dots f_m$ equal to zero. This is in contrast to the case of regular tilings by $\mathcal{T}_n, n$ even, discussed in~\cite{nitica-L-shaped}, which is very well understood.

For completeness, we also consider the tiling set $\tilde {\mathcal{T}}_n$, consisting of $\mathcal{T}_n, n$ odd, and an extra $2\times 2$ square.

\begin{thm}\label{thm:last63} If $n\ge 5$ odd, any region in a square lattice can be signed tiled by $\tilde {\mathcal{T}}_n$.
\end{thm}

Theorem~\ref{thm:last63} is proved in Section 7.

Barnes developed in~\cite{Barnes1, Barnes2} a general method for solving signed tiling problems with complex weights. In Section 8 we review the method of Barnes and offer an alternative proof of Theorem~\ref{thm-main} based on this method. Having available a Gr\"obner basis for the tiling set helps even if Barnes method is used.

\begin{thm}\label{thm-main-barnes} If complex or rational weights are allowed to replace the integral weights $\pm 1$, a rectangle can be signed tiled by $\mathcal{T}_n, n\ge 5$ odd, if and only if it has a side divisible by $n$.
\end{thm}

Signed tilings by $\mathcal{T}_n, n$ even, are more complicated then in the odd case. We will discussed them in a forthcoming paper~\cite{gill-nitica-even}.

A final comment about the paper. While the methods that we use are well known, and algorithmic when applied to a particular tiling problem, here we apply them uniformly to solve simultaneously an infinite collection of problems.

\section{Summary of Gr\"obner basis theory}

An introduction to signed tilings can be found in the paper of Conway and Lagarias~\cite{con-lag}. One investigates there signed tilings by the $3$-bone, a tile consisting of three adjacent regular hexagons. The Gr\"obner basis approach to signed polyomino tilings was proposed by Bodini and Nouvel~\cite{B-N}. In~\cite{n-bone} one uses this approach to study signed tilings by the $n$-bone.

Let $R[\underline{X}] = R[X_1, \dots , X_k]$ be the ring of polynomials with coefficients in a principal ideal domain (PID) $R$. The only (PID) of interest in this paper is $\mathbb{Z}$, the ring of integers. A \emph{term} in the variables $x_1,\dots,x_k$ is a power product $x_1^{\alpha_1}x_2^{\alpha_2}\dots x_{\ell}^{\alpha_{\ell}}$ with $\alpha_i\in \mathbb{N}, 1\le i\le \ell$; in particular
$1=x_1^0\dots x_{\ell}^0$ is a term. A term with an associated coefficient from $R$ is called \emph{monomial}. We endow the set of terms with the \emph{total degree-lexicographical order}, in which we first compare the degrees of the monomials and then break the ties by means of lexicographic order for the order $x_1>x_2>\dots>x_{\ell}$ on the variables. If the variables are only $x,y$ and $x>y,$ this gives the total order:
\begin{equation}
1<y<x<y^2<xy<x^2<y^3<xy^2<x^2y<x^3<y^4<\cdots.
\end{equation}
For $P\in R[\underline{X}]$ we denote by $HT(P)$ the leading term in $P$ with respect to the above order and by $HM(P)$ the monomial of $HT(P)$. We denote by $HC(P)$ the coefficient of the leading monomial in $P$. We denote by $T(P)$ the set of terms appearing in $P$, which we assume to be in simplest form. We denote by $M(P)$ the set of monomials in $P$. For a given ideal $I\subseteq R[\underline{X}]$ an associated Gr\"obner basis may be introduced for example as in \cite[Chapters 5, 10]{B}. Our summary follows the approach there. If $G\subseteq R[\underline{X}]$ is a finite set, we denote by $I(G)$ the ideal generated by $G$ in $R[\underline{X}]$.

\begin{definition} Let $f,g,p\in R[\underline{X}]$. We say that $f$ $D$-reduces to $g$ modulo $p$ and write $f \underset{p}{\to} g$ if there exists $m\in M(f)$ with $HM(p)\vert m$, say $m=m'\cdot HM(p),$ and $g=f-m'p$. For a finite set $G\subseteq R[\underline{X}]$, we denote by $\overset{*}{\underset{G}{\to}}$ the reflexive-transitive closure of $\underset{p}{\to}, p\in G$. We say that $g$ is a normal form for $f$ with respect to $G$ if $f \overset{*}{\underset{G}{\to}} g$ and no further $D$-reduction is possible. We say that $f$ is $D$-reducible modulo $G$ if $f\overset{*}{\underset{G}{\to}} 0$.
\end{definition}

It is clear that if $f\overset{*}{\underset{G}{\to}} 0$, then $f$ belongs to the ideal generated by $G$ in $R[\underline{X}]$. The converse is also true if $G$ is a Gr\"obner basis.

\begin{definition} A $D$-Gr\"obner basis is a finite set $G$ of $R[\underline{X}]$ with the property that all $D$-normal forms modulo $G$ of elements of $I(G)$ equal zero. If $I\subseteq R[\underline{X}]$ is an ideal, then a $D$-Gr\"oebner basis of $I$ is a $D$-Gr\"oebner basis that generates the ideal $I$.
\end{definition}

\begin{propo} Let $G$ be a finite set of $R[\underline{X}]$. Then the following statements are equivalent:
\begin{enumerate}
\item $G$ is a Gr\"oebner basis.
\item Every $f\not = 0, f\in I(G),$ is $D$-reducible modulo $G$.
\end{enumerate}
\end{propo}

Note that if $R$ is only a (PID), the normal form of the division of $f$ by $G$ is not unique.

We introduce now the notions of $S$-polynomial and $G$-polynomial.

\begin{definition} Let $0\not = g_i\in R[\underline{X}], i=1,2,$ with $HC(g_i)=a_i$ and $HT(g_i)=t_i$. Let $a=b_ia_i=\text{lcm}(a_1,a_2)$ with $b_i\in R$, and $t=s_it_i=\text{lcm}(t_1,t_2)$ with $s_i\in T$. Let $c_1,c_2\in R$ such that $\text{gcd}(a_1,a_2)=c_1a_1+c_2a_2$. Then:
\begin{equation}
\begin{gathered}
S(g_1,g_2)=b_1s_1g_1-b_2s_2g_2,\\
G(g_1,g_2)=c_1s_1g_1+c_2s_2g_2.
\end{gathered}
\end{equation}
\end{definition}

\begin{rmk} If $HC(g_1)=HC(g_2)$, then $G(g_1,g_2)$ can be chosen to be $g_1$.
\end{rmk}

\begin{thm} Let $G$ be a finite set of $R[\underline{X}]$. Assume that for all $g_1,g_2\in G$, $S(g_1,g_2)\overset{*}{\underset{G}{\to}} 0$ and $G(g_1,g_2)$ is top-$D$-reducible modulo $G$. Then $G$ is a Gr\"obner basis.
\end{thm}

Assume now that $R$ is an Euclidean domain with unique remainders (see~\cite[p. 463]{B}). This is the case for the ring of integers $\mathbb{Z}$ if we specify remainders upon division by $0\not =m$ to be in the interval $[0,m)$.

\begin{definition} Let $f,g,p\in R[\underline{X}]$. We say that $f$ $E$-reduces to $g$ modulo $p$ and write $f \underset{p}{\to} g$ if there exists $m=at\in M(f)$ with $HM(p)\vert t$, say $t=s\cdot HT(p),$ and $g=f-qsp$ where $0\not =q\in R$ is the quotient of $a$ upon division with unique remainder by $HC(p)$.
\end{definition}

\begin{propo} $E$-reduction extends $D$-reduction, i.e., every $D$-reduction step in an $E$-reduction step.
\end{propo}

\begin{thm} Let $R$ be an Euclidean domain with unique remainders, and assume $G\subseteq R[\underline{X}]$ is a $D$-Gr\"obner basis. Then the following hold:
\begin{enumerate}
\item $f \overset{*}{\underset{G}{\to}} 0$ for all $f\in I(G),$ where $\overset{*}{\underset{G}{\to}}$ denotes the $E$-reduction modulo $G$.
\item $E$-reduction modulo $G$ has unique normal forms.
\end{enumerate}
\end{thm}

The following result connect signed tilings and Gr\"obner bases. See~\cite{B-N} and~\cite{n-bone} for a proof.

\begin{thm}\label{thm-groebner-tiling} A polyomino $P$ admits a signed tiling by translates of prototiles
$P_1, P_2, \dots , P_k$ if and only if for some (test) monomial $x^{\alpha}y^{\beta}$ the polynomial
$x^{\alpha}y^{\beta}f_P$ is in the ideal generated in $\mathbb{Z}[X,Y]$ by the polynomials $f_{P_1},\dots, f_{P_k}$.
Moreover, the set of test monomials $\mathcal{T} = \{x^{\alpha}\}$ can be chosen from any set
$T\subseteq\mathbb{N}^n$ of multi-indices which is cofinal in $(\mathbb{N}^n,\le)$.
\end{thm}

\section{Gr\"obner basis for $\mathcal{T}_n, n$ odd}

We write $n=2k+1,$ $k\ge 2$. The polynomials (in a condensed form) associated to the tiles in $\mathcal{T}_n$ are:

\begin{equation}\label{eq:generators-k}
\begin{aligned}
G_1(k)&=\frac{y^{2k}-1}{y-1}+x,
G_2(k)=y^{2k-1}+x\cdot\frac{y^{2k}-1}{y-1},
G_3(k)=y+\frac{x^{2k}-1}{x-1},
G_4(k)=y\cdot\frac{x^{2k}-1}{x-1}+x^{2k-1}.
\end{aligned}
\end{equation}

We show in the rest of this section that a Gr\"oebner basis for the ideal generated in $\mathbb{Z}[X,Y]$ by $G_1(k),$ $G_2(k),$ $G_3(k),G_4(k),$ is given by the polynomials (written in condensed from):

\begin{equation}\label{eq:basis-k}
\begin{aligned}
B_1(k)&=\frac{y^{k+2}-1}{y-1}+x\cdot\frac{x^{k-1}-1}{x-1}, \
B_2(k)=\frac{y^{k+1}-1}{y-1}+x\cdot\frac{x^{k}-1}{x-1}. \
B_3(k)=xy-1.
\end{aligned}
\end{equation}

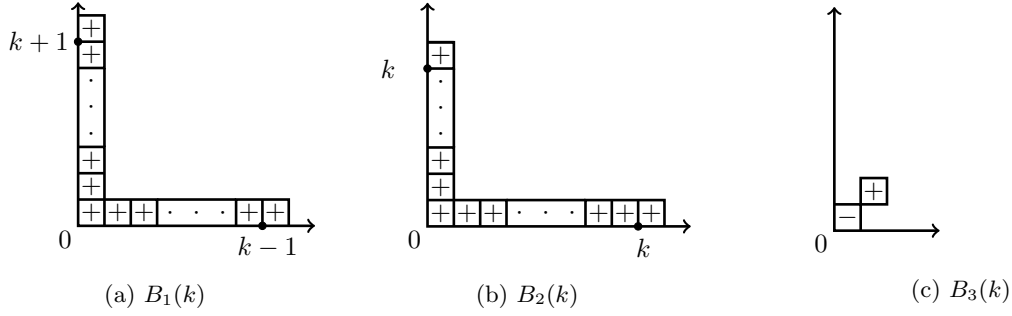
\begin{figure}[h!]
\centering
\begin{subfigure}{.25\textwidth}
\begin{tikzpicture}[scale=.35]
\draw [line width = 1,<->] (0,8.5)--(0,0)--(9,0);
\draw [line width = 1] (0,1)--(3,1)--(3,0);
\draw [line width = 1] (1,0)--(1,3)--(0,3);
\draw [line width = 1] (2,0)--(2,1);
\draw [line width = 1] (0,2)--(1,2);
\draw [line width = 1] (6,0)--(6,1)--(8,1)--(8,0);
\draw [line width = 1] (0,6)--(1,6)--(1,8)--(0,8);
\draw [line width = 1] (7,0)--(7,1);
\draw [line width = 1] (0,7)--(1,7);
\draw [line width = 1] (1,3)--(1,6);
\draw [line width = 1] (3,1)--(6,1);

\node at (.5,.5) {$+$};
\node at (1.5,.5) {$+$};
\node at (2.5,.5) {$+$};
\node at (6.5,.5) {$+$};
\node at (7.5,.5) {$+$};

\node at (.5,1.5) {$+$};
\node at (.5,2.5) {$+$};
\node at (.5,6.5) {$+$};
\node at (.5,7.5) {$+$};

\node at (3.5,.5) {$\cdot$};
\node at (4.5,.5) {$\cdot$};
\node at (5.5,.5) {$\cdot$};

\node at (.5,3.5) {$\cdot$};
\node at (.5,4.5) {$\cdot$};
\node at (.5,5.5) {$\cdot$};

\draw [fill] (7,0) circle (4pt);
\draw [fill] (0,7) circle (4pt);

\node at (-.5, -.5) {$0$};
\node at (7.2, -.8) {$k-1$};
\node at (-1.5, 7) {$k+1$};
\end{tikzpicture}
\caption{$B_1(k)$}
\end{subfigure}
~~~~~
\begin{subfigure}{.25\textwidth}
\begin{tikzpicture}[scale=.35]
\draw [line width = 1,<->] (0,8.5)--(0,0)--(10,0);
\draw [line width = 1] (0,1)--(3,1)--(3,0);
\draw [line width = 1] (1,0)--(1,3)--(0,3);
\draw [line width = 1] (2,0)--(2,1);
\draw [line width = 1] (0,2)--(1,2);
\draw [line width = 1] (6,0)--(6,1)--(9,1)--(9,0);
\draw [line width = 1] (0,6)--(1,6)--(1,7)--(0,7);
\draw [line width = 1] (7,0)--(7,1);
\draw [line width = 1] (8,0)--(8,1);
\draw [line width = 1] (0,7)--(1,7);
\draw [line width = 1] (1,3)--(1,6);
\draw [line width = 1] (3,1)--(6,1);

\node at (.5,.5) {$+$};
\node at (1.5,.5) {$+$};
\node at (2.5,.5) {$+$};
\node at (6.5,.5) {$+$};
\node at (7.5,.5) {$+$};

\node at (.5,1.5) {$+$};
\node at (.5,2.5) {$+$};
\node at (.5,6.5) {$+$};

\node at (8.5,.5) {$+$};

\node at (3.5,.5) {$\cdot$};
\node at (4.5,.5) {$\cdot$};
\node at (5.5,.5) {$\cdot$};

\node at (.5,3.5) {$\cdot$};
\node at (.5,4.5) {$\cdot$};
\node at (.5,5.5) {$\cdot$};

\draw [fill] (8,0) circle (4pt);
\draw [fill] (0,6) circle (4pt);

\node at (-.5, -.5) {$0$};
\node at (8.2, -.8) {$k$};
\node at (-1.5, 6) {$k$};
\end{tikzpicture}
\caption{$B_2(k)$}
\end{subfigure}
~~~~~~~~~~~~
\begin{subfigure}{.25\textwidth}
\begin{tikzpicture}[scale=.35]
\draw [line width = 1,<->] (0,8.5)--(0,0)--(4,0);
\draw [line width = 1] (0,1)--(1,1)--(1,0);
\draw [line width = 1] (1,1)--(1,2)--(2,2)--(2,1)--(1,1);

\node at (.5,.5) {$-$};
\node at (1.5,1.5) {$+$};
\node at (-.5, -.5) {$0$};
\end{tikzpicture}
\caption{$B_3(k)$}
\end{subfigure}
\caption{The Gr\"obner basis $\{B_1(k), B_2(k), B_3(k)\}.$}
\label{fig:signed-tiles-odd}
\end{figure}

It is convenient for us to look at the elements of the basis geometrically, as signed tiles, see Figure~\ref{fig:signed-tiles-odd}.

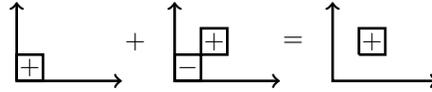
\begin{figure}[h!]
\centering
\begin{tikzpicture}[scale=.35]
\draw [line width = 1,<->] (0,3)--(0,0)--(4,0);
\draw [line width = 1] (0,0)--(0,1)--(1,1)--(1,0)--(0,0);

\node at (.5, .5) {$+$};

\node at (4.5, 1.5) {$+$};

\draw [line width = 1,<->] (6,3)--(6,0)--(10,0);
\draw [line width = 1] (6,0)--(7,0)--(7,2)--(8,2)--(8,1)--(6,1)--(6,0);

\node at (6.5, .5) {$-$};
\node at (7.5, 1.5) {$+$};

\node at (10.5, 1.5) {$=$};

\draw [line width = 1,<->] (12,3)--(12,0)--(16,0);
\draw [line width = 1] (13,1)--(14,1)--(14,2)--(13,2)--(13,1);

\node at (13.5, 1.5) {$+$};

\end{tikzpicture}
\caption{Tiles arithmetic.}
\label{fig:tilesarth3}
\end{figure}

The presence of $B_3(k)$ in the basis allows to reduce the algebraic proofs to combinatorial considerations. Indeed, using addition by a multiple of $B_3(k)$, one can translate, along a vector parallel to the first bisector $y=x$, cells labeled by $+1$ from one position in the square lattice to another. See Figure~\ref{fig:tilesarth3}.

We will use this property repeatedly to check certain algebraic identities.

\begin{propo}\label{p:prop5} $G_1(k), G_2(k), G_3(k), G_4(k)$ are in the ideal generated by $B_1(k),$ $B_2(k),$ $B_3(k)$.
\end{propo}

\begin{pf} The geometric proofs appear in Figure~\ref{fig:proofs1}. First we translate one of the tiles $B_i(k), i=1,2,$ multiplying by a power of $x$ or a power of $y$, and then rearrange the cells from $B_i(k)$ using diagonal translations given by multiples of $B_3(k)$. The initial tiles $B_i(k), i=1,2,$ have the cells marked by a $+$, and the final tiles $G_i(k), i=1,2,3,4,$ are colored in light gray.

\begin{figure}[h!]
\centering
\begin{subfigure}{.25\textwidth}
\begin{tikzpicture}[scale=.33]
\draw [line width = 1, fill=lightgray] (0,6)--(0,14)--(1,14)--(1,6)--(0,6);
\draw [line width = 1,<->] (0,15)--(0,0)--(9,0);
\draw [line width = 1] (0,7)--(3,7)--(3,6);
\draw [line width = 1] (1,6)--(1,9)--(0,9);
\draw [line width = 1] (2,6)--(2,7);
\draw [line width = 1] (0,8)--(1,8);
\draw [line width = 1] (6,6)--(6,7)--(8,7)--(8,6);
\draw [line width = 1] (0,12)--(1,12)--(1,14)--(0,14);
\draw [line width = 1] (7,6)--(7,7);
\draw [line width = 1] (0,13)--(1,13);
\draw [line width = 1] (0,6)--(3,6);
\draw [line width = 1] (6,6)--(8,6);
\draw [line width = 1] (1,9)--(1,12);
\draw [line width = 1] (3,7)--(6,7);
\draw [line width = 1] (3,6)--(6,6);

\draw [line width = 1, fill=lightgray] (0,6)--(1,6)--(1,1)--(2,1)--(2,0)--(0,0)--(0,6);
\draw [line width = 1] (1,0)--(1,1);
\draw [line width = 1] (0,1)--(1,1);
\draw [line width = 1] (0,4)--(1,4);
\draw [line width = 1] (0,5)--(1,5);

\draw [line width = 1, ->] (1.5,6.5)--(.5,5.5);
\draw [line width = 1,->] (2.5,6.5)--(.5,4.5);
\draw [line width = 1, ->] (6.5,6.5)--(.5,.5);
\draw [line width = 1,->] (7.5,6.5)--(1.5,.5);

\node at (.5,6.5) {$+$};
\node at (1.5,6.5) {$+$};
\node at (2.5,6.5) {$+$};
\node at (6.5,6.5) {$+$};
\node at (7.5,6.5) {$+$};

\node at (.5,7.5) {$+$};
\node at (.5,8.5) {$+$};
\node at (.5,12.5) {$+$};
\node at (.5,13.5) {$+$};

\node at (3.5,6.5) {$\cdot$};
\node at (4.5,6.5) {$\cdot$};
\node at (5.5,6.5) {$\cdot$};

\node at (.5,9.5) {$\cdot$};
\node at (.5,10.5) {$\cdot$};
\node at (.5,11.5) {$\cdot$};

\node at (.5,3.5) {$\cdot$};
\node at (.5,2.5) {$\cdot$};
\node at (.5,1.5) {$\cdot$};

\draw [fill] (7,0) circle (4pt);
\draw [fill] (0,6) circle (4pt);
\draw [fill] (0,13) circle (4pt);

\node at (-.5, -.5) {$0$};
\node at (7.2, -.8) {$k-1$};
\node at (-1.5, 6) {$k-2$};
\node at (-1.5, 13) {$2k-1$};
\end{tikzpicture}
\caption{$G_1(k)$}
\end{subfigure}
~~~~
\begin{subfigure}{.25\textwidth}
\begin{tikzpicture}[scale=.33]
\draw [line width = 1, fill=lightgray] (1,0)--(1,13)--(0,13)--(0,14)--(2,14)--(2,0);
\draw [line width = 1,<->] (0,15)--(0,0)--(9,0);
\draw [line width = 1] (0,7)--(3,7)--(3,6);
\draw [line width = 1] (1,6)--(1,9)--(0,9);
\draw [line width = 1] (2,6)--(2,7);
\draw [line width = 1] (0,8)--(1,8);
\draw [line width = 1] (6,6)--(6,7)--(8,7)--(8,6);
\draw [line width = 1] (0,12)--(1,12)--(1,14)--(0,14);
\draw [line width = 1] (7,6)--(7,7);
\draw [line width = 1] (0,13)--(1,13);
\draw [line width = 1] (0,6)--(3,6);
\draw [line width = 1] (6,6)--(8,6);
\draw [line width = 1] (1,9)--(1,12);
\draw [line width = 1] (3,7)--(6,7);
\draw [line width = 1] (3,6)--(6,6);

\draw [line width = 1] (1,1)--(2,1);
\draw [line width = 1] (1,2)--(2,2);
\draw [line width = 1] (1,8)--(2,8);
\draw [line width = 1] (1,5)--(2,5);
\draw [line width = 1] (1,9)--(2,9);
\draw [line width = 1] (1,10)--(2,10);
\draw [line width = 1] (1,13)--(2,13);

\draw [line width = 1, <-] (1.5,7.5)--(.5,6.5);
\draw [line width = 1, <-] (1.5,8.5)--(.5,7.5);
\draw [line width = 1, <-] (1.5,9.5)--(.5,8.5);
\draw [line width = 1, <-] (1.5,13.5)--(.5,12.5);
\draw [line width = 1,->] (2.5,6.5)--(1.5,5.5);
\draw [line width = 1, ->] (6.5,6.5)--(1.5,1.5);
\draw [line width = 1,->] (7.5,6.5)--(1.5,.5);

\node at (.5,6.5) {$+$};
\node at (1.5,6.5) {$+$};
\node at (2.5,6.5) {$+$};
\node at (6.5,6.5) {$+$};
\node at (7.5,6.5) {$+$};

\node at (.5,7.5) {$+$};
\node at (.5,8.5) {$+$};
\node at (.5,12.5) {$+$};
\node at (.5,13.5) {$+$};

\node at (3.5,6.5) {$\cdot$};
\node at (4.5,6.5) {$\cdot$};
\node at (5.5,6.5) {$\cdot$};

\node at (.5,9.5) {$\cdot$};
\node at (.5,10.5) {$\cdot$};
\node at (.5,11.5) {$\cdot$};

\node at (1.5,12.5) {$\cdot$};
\node at (1.5,10.5) {$\cdot$};
\node at (1.5,11.5) {$\cdot$};

\node at (1.5,3.5) {$\cdot$};
\node at (1.5,2.5) {$\cdot$};
\node at (1.5,4.5) {$\cdot$};

\draw [fill] (7,0) circle (4pt);
\draw [fill] (0,6) circle (4pt);
\draw [fill] (0,13) circle (4pt);

\node at (-.5, -.5) {$0$};
\node at (7.2, -.8) {$k-1$};
\node at (-1.5, 6) {$k-1$};
\node at (-1.5, 13) {$2k-1$};
\end{tikzpicture}
\caption{$G_2(k)$}
\end{subfigure}
~~~~
\begin{subfigure}{.25\textwidth}
\begin{tikzpicture}[scale=.33]
\draw [line width = 1, fill=lightgray] (0,2)--(1,2)--(1,1)--(14,1)--(14,0)--(0,0)--(0,2);

\draw [line width = 1,<->] (0,14)--(0,0)--(15,0);
\draw [line width = 1] (5,1)--(8,1)--(8,0);
\draw [line width = 1] (6,0)--(6,3)--(5,3);
\draw [line width = 1] (7,0)--(7,1);
\draw [line width = 1] (5,2)--(6,2);
\draw [line width = 1] (11,0)--(11,1)--(14,1)--(14,0);
\draw [line width = 1] (5,6)--(6,6)--(6,7)--(5,7);
\draw [line width = 1] (12,0)--(12,1);
\draw [line width = 1] (13,0)--(13,1);
\draw [line width = 1] (5,7)--(6,7);
\draw [line width = 1] (5,0)--(5,7);
\draw [line width = 1] (6,3)--(6,6);
\draw [line width = 1] (0,1)--(1,1);
\draw [line width = 1] (3,0)--(3,1);
\draw [line width = 1] (4,0)--(4,1);
\draw [line width = 1,->] (5.5,1.5)--(4.5,.5);
\draw [line width = 1,->] (5.5,2.5)--(3.5,.5);
\draw [line width = 1,->] (5.5,6.5)--(.5,1.5);

\node at (5.5,.5) {$+$};
\node at (6.5,.5) {$+$};
\node at (7.5,.5) {$+$};
\node at (11.5,.5) {$+$};
\node at (12.5,.5) {$+$};

\node at (5.5,1.5) {$+$};
\node at (5.5,2.5) {$+$};
\node at (5.5,6.5) {$+$};

\node at (13.5,.5) {$+$};

\node at (8.5,.5) {$\cdot$};
\node at (9.5,.5) {$\cdot$};
\node at (10.5,.5) {$\cdot$};

\node at (.5,.5) {$\cdot$};
\node at (1.5,.5) {$\cdot$};
\node at (2.5,.5) {$\cdot$};

\node at (5.5,3.5) {$\cdot$};
\node at (5.5,4.5) {$\cdot$};
\node at (5.5,5.5) {$\cdot$};

\draw [fill] (5,0) circle (4pt);
\draw [fill] (0,6) circle (4pt);
\draw [fill] (13,0) circle (4pt);

\node at (-.5, -.5) {$0$};
\node at (5, -.8) {$k-1$};
\node at (-1.5, 6) {$k$};
\node at (13, -.8) {$2k-1$};
\end{tikzpicture}
\caption{$G_3(k)$}
\end{subfigure}
~~~~
\begin{subfigure}{.25\textwidth}
\begin{tikzpicture}[scale=.33]
\draw [line width = 1, fill=lightgray] (0,2)--(14,2)--(14,0)--(13,0)--(13,1)--(0,1)--(0,2);

\draw [line width = 1,<->] (0,8)--(0,0)--(15,0);
\draw [line width = 1] (5,1)--(8,1)--(8,0);
\draw [line width = 1] (6,0)--(6,3)--(5,3);
\draw [line width = 1] (7,0)--(7,1);
\draw [line width = 1] (5,2)--(6,2);
\draw [line width = 1] (11,0)--(11,1)--(14,1)--(14,0);
\draw [line width = 1] (5,6)--(6,6)--(6,7)--(5,7);
\draw [line width = 1] (12,0)--(12,1);
\draw [line width = 1] (13,0)--(13,1);
\draw [line width = 1] (5,7)--(6,7);
\draw [line width = 1] (5,0)--(5,7);
\draw [line width = 1] (6,3)--(6,6);
\draw [line width = 1] (0,1)--(1,1);
\draw [line width = 1] (1,2)--(1,1);
\draw [line width = 1] (4,2)--(4,1);
\draw [line width = 1] (7,2)--(7,1);
\draw [line width = 1] (8,2)--(8,1);
\draw [line width = 1] (9,2)--(9,1);
\draw [line width = 1] (12,2)--(12,1);
\draw [line width = 1] (13,2)--(13,1);
\draw [line width = 1,<-] (6.5,1.5)--(5.5,.5);
\draw [line width = 1,<-] (7.5,1.5)--(6.5,.5);
\draw [line width = 1,<-] (8.5,1.5)--(7.5,.5);
\draw [line width = 1,->] (5.5,2.5)--(4.5,1.5);
\draw [line width = 1,->] (5.5,6.5)--(.5,1.5);
\draw [line width = 1,<-] (12.5,1.5)--(11.5,.5);
\draw [line width = 1,<-] (13.5,1.5)--(12.5,.5);

\node at (5.5,.5) {$+$};
\node at (6.5,.5) {$+$};
\node at (7.5,.5) {$+$};
\node at (11.5,.5) {$+$};
\node at (12.5,.5) {$+$};

\node at (5.5,1.5) {$+$};
\node at (5.5,2.5) {$+$};
\node at (5.5,6.5) {$+$};

\node at (13.5,.5) {$+$};

\node at (11.5,1.5) {$\cdot$};
\node at (9.5,1.5) {$\cdot$};
\node at (10.5,1.5) {$\cdot$};

\node at (8.5,.5) {$\cdot$};
\node at (9.5,.5) {$\cdot$};
\node at (10.5,.5) {$\cdot$};

\node at (1.5,1.5) {$\cdot$};
\node at (2.5,1.5) {$\cdot$};
\node at (3.5,1.5) {$\cdot$};

\node at (5.5,3.5) {$\cdot$};
\node at (5.5,4.5) {$\cdot$};
\node at (5.5,5.5) {$\cdot$};

\draw [fill] (5,0) circle (4pt);
\draw [fill] (0,6) circle (4pt);
\draw [fill] (13,0) circle (4pt);

\node at (-.5, -.5) {$0$};
\node at (5, -.8) {$k$};
\node at (-1.5, 6) {$k$};
\node at (13, -.8) {$2k-1$};
\end{tikzpicture}
\caption{$G_4(k)$}
\end{subfigure}
\caption{Generating $\{G_1(k),G_2(k),G_3(k),G_4(k)\}$ from $\{B_1(k), B_2(k), B_3(k)\}.$}
\label{fig:proofs1}
\end{figure}
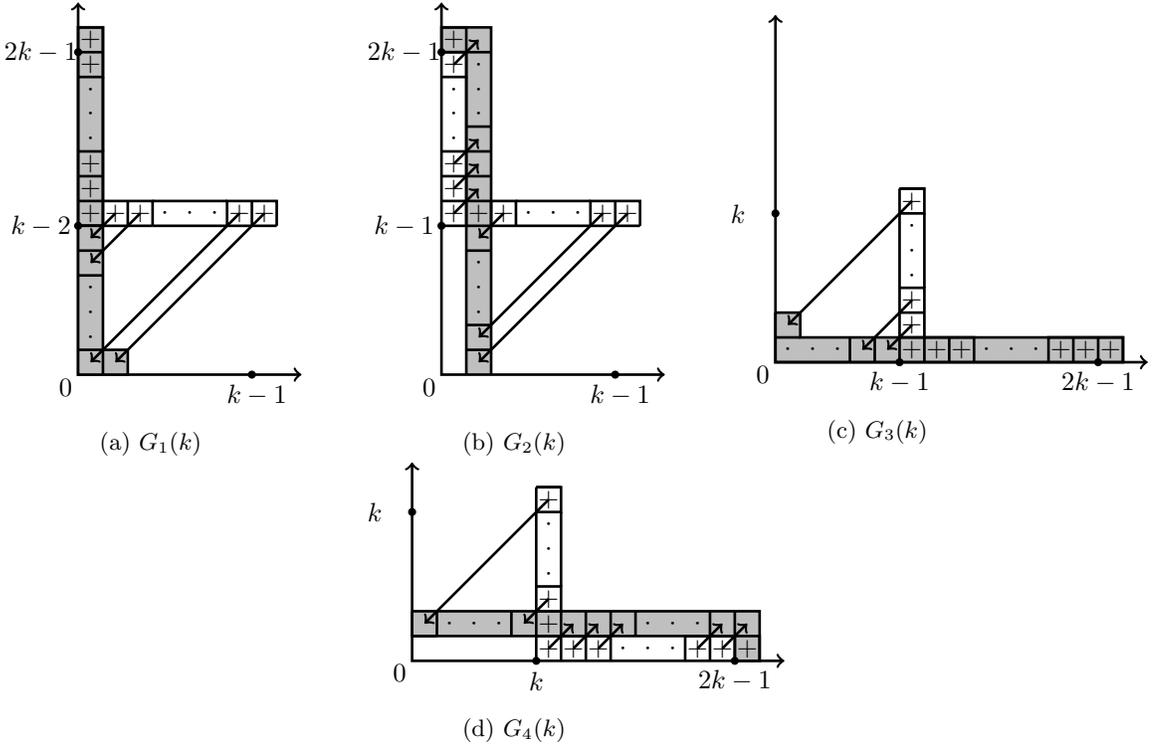
\end{pf}

\begin{propo}\label{p:prop6} $B_1(k), B_2(k), B_3(k)$ are in the ideal generated by $G_1(k),$ $G_2(k),$ $G_3(k),$ $G_4(k)$.
\end{propo}

\begin{pf} We first show that $B_3$ belongs to the ideal generated by $G_1(k), G_2(k), G_3(k), G_4(k)$. One has:
\begin{eqnarray}\label{eq:eq1207}
\begin{aligned}
B_3(k)&=-G_1(k)+G_2(k)+\left ( -xy\cdot \frac{y^{2(k-1)}-1}{y-1}+y\cdot \frac{y^{2(k-1)}-1}{y^2-1} \right )G_3(k)+xy\cdot \frac{y^{2(k-1)}-1}{y^2-1}\cdot G_4(k).
\end{aligned}
\end{eqnarray}

Using~\eqref{eq:generators-k}, the RHS of equation~\eqref{eq:eq1207} becomes:
\begin{equation*}
\begin{aligned}
&\frac{1}{(y^2-1)(x-1)}\Big [ -(y^{2k}-1)(y+1)(x-1)-x(y^2-1)(x-1)+y^{2k-1}(y^2-1)(x-1)\\
&+x(y^{2k}-1)(y+1)(x-1)+ [-xy(y+1)(y^{2(k-1)}-1)+y(y^{2(k-1)}-1)]\cdot [y(x-1)+x^{2k}-1]\\
&+xy(y^{2(k-1)}-1)\cdot [y(x^{2k}-1)+x^{2k-1}(x-1)]\Big ]\\
=&\frac{1}{(y^2-1)(x-1)}\Big [ (-y^{2k}+1)(xy-y+x-1)-x(y^2x-y^2-x+1)+y^{2k-1}(y^2x-y^2-x+1)\\
&+(xy^{2k}-x)(xy-y+x-1)+(-xy^{2k}+xy^2-xy^{2k-1}+xy+y^{2k-1}-y)\\
&\cdot (xy-y+x^{2k}-1)+(xy^{2k-1}-xy)(yx^{2k}-y+x^{2k}-x^{2k-1}) \Big ]
\end{aligned}
\end{equation*}

\begin{equation*}
\begin{aligned}
=&\frac{1}{(y^2-1)(x-1)}\Big [ -y^{2k+1}x+y^{2k+1}-y^{2k}x+y^{2k}+xy-y+x-1-x^2y^2+xy^2+x^2-x\\
&+xy^{2k+1}-y^{2k+1}-xy^{2k-1}+y^{2k-1}+x^2y^{2k+1}-xy^{2k+1}+x^2y^{2k}-xy^{2k}\\
&-x^2y+xy-x^2+x-x^2y^{2k+1}+xy^{2k+1}-x^{2k+1}y^{2k}+xy^{2k}\\
&+x^2y^3-xy^3+x^{2k+1}y^2-xy^2-x^2y^{2k}+xy^{2k}-x^{2k+1}y^{2k-1}+xy^{2k-1}\\
&+x^2y^2-xy^2+x^{2k+1}y-xy+xy^{2k}-y^{2k}+x^{2k}y^{2k-1}-y^{2k-1}\\
&-xy^2+y^2-x^{2k}y+y+x^{2k+1}y^{2k}-x^{2k+1}y^2-xy^{2k}+xy^2+x^{2k+1}y^{2k-1}-x^{2k+1}y\\
&-x^{2k}y^{2k-1}+x^{2k}y\Big ]\\
=&\frac{1}{(y^2-1)(x-1)}\Big [xy-1-x^2y+x^2y^3-xy^3-xy^2+y^2+x\Big ]=xy-1=B_3(k).
\end{aligned}
\end{equation*}

After we obtain $B_3(k)$, polynomials $B_1(k), B_2(k)$ can be obtained geometrically by reversing the processes in Figure~\ref{fig:proofs1}. Reversing the process in Figure~\ref{fig:proofs1}, a), we first obtain a copy of $y^{k-2}B_1(k)$. This copy can be translated to the right using multiplication by $x^{k-2}$, and then can be pulled back with the corner in the origin using a translation by a vector parallel to $y=x$. Reversing the process in Figure~\ref{fig:proofs1}, c), we first obtain a copy of $x^{k-1}B_2(k)$. This copy can be translated up uisng multiplication by $y^{k-1}$, and then can be pulled back with the corner in the origin using a translation by a vector parallel to $y=x$.

A step by step geometric proof of formula~\ref{eq:eq1207} for $n=7$ is shown in Figure~\ref{fig:special-form-33}. All cells in the square lattice without any label have weight zero. The proof can be easily generalized for any odd $n$.

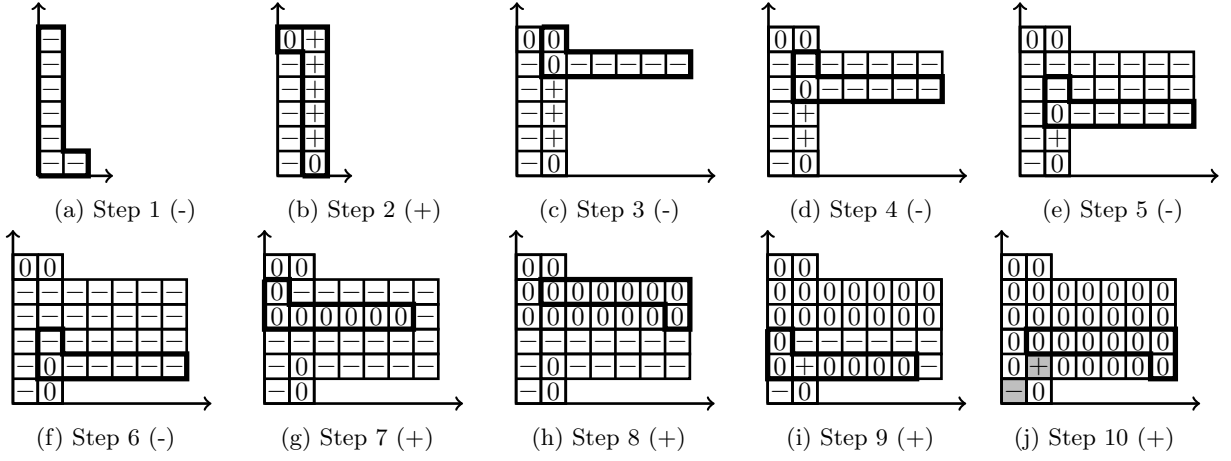
\begin{figure}[h!]
\centering
\begin{subfigure}{.15\textwidth}
\begin{tikzpicture}[scale=.33]

\draw [line width = 1,<->] (0,7)--(0,0)--(3,0);
\draw [line width = 2] (0,0)--(2,0)--(2,1)--(1,1)--(1,6)--(0,6)--(0,0);
\draw [line width = 1] (1,0)--(1,1);
\draw [line width = 1] (0,1)--(1,1);
\draw [line width = 1] (0,2)--(1,2);
\draw [line width = 1] (0,3)--(1,3);
\draw [line width = 1] (0,4)--(1,4);
\draw [line width = 1] (0,5)--(1,5);

\node at (.5,.5) {$-$};
\node at (.5,1.5) {$-$};
\node at (.5,2.5) {$-$};
\node at (.5,3.5) {$-$};
\node at (.5,4.5) {$-$};
\node at (.5,5.5) {$-$};

\node at (1.5,.5) {$-$};
\end{tikzpicture}
\caption{Step 1 (-)}
\end{subfigure}
~~~~
\begin{subfigure}{.15\textwidth}
\begin{tikzpicture}[scale=.33]

\draw [line width = 1,<->] (0,7)--(0,0)--(3,0);
\draw [line width = 2] (1,0)--(2,0)--(2,6)--(0,6)--(0,5)--(1,5)--(1,0);
\draw [line width = 1] (0,0)--(2,0)--(2,1)--(1,1)--(1,6)--(0,6)--(0,0);
\draw [line width = 1] (1,0)--(1,1);
\draw [line width = 1] (0,1)--(1,1);
\draw [line width = 1] (0,2)--(1,2);
\draw [line width = 1] (0,3)--(1,3);
\draw [line width = 1] (0,4)--(1,4);
\draw [line width = 1] (0,5)--(1,5);

\draw [line width = 1] (1,2)--(2,2);
\draw [line width = 1] (1,3)--(2,3);
\draw [line width = 1] (1,4)--(2,4);
\draw [line width = 1] (1,5)--(2,5);

\node at (1.5,1.5) {$+$};
\node at (1.5,2.5) {$+$};
\node at (1.5,3.5) {$+$};
\node at (1.5,4.5) {$+$};
\node at (1.5,5.5) {$+$};

\node at (.5,.5) {$-$};
\node at (.5,1.5) {$-$};
\node at (.5,2.5) {$-$};
\node at (.5,3.5) {$-$};
\node at (.5,4.5) {$-$};
\node at (.5,5.5) {$0$};

\node at (1.5,.5) {$0$};
\end{tikzpicture}
\caption{Step 2 (+)}
\end{subfigure}
~~~~
\begin{subfigure}{.16\textwidth}
\begin{tikzpicture}[scale=.33]

\draw [line width = 1,<->] (0,7)--(0,0)--(8,0);
\draw [line width = 2] (1,6)--(2,6)--(2,5)--(7,5)--(7,4)--(1,4)--(1,6);
\draw [line width = 1] (1,0)--(2,0)--(2,6)--(0,6)--(0,5)--(1,5)--(1,0);
\draw [line width = 1] (0,0)--(2,0)--(2,1)--(1,1)--(1,6)--(0,6)--(0,0);
\draw [line width = 1] (1,0)--(1,1);
\draw [line width = 1] (0,1)--(1,1);
\draw [line width = 1] (0,2)--(1,2);
\draw [line width = 1] (0,3)--(1,3);
\draw [line width = 1] (0,4)--(1,4);
\draw [line width = 1] (0,5)--(1,5);

\draw [line width = 1] (1,2)--(2,2);
\draw [line width = 1] (1,3)--(2,3);
\draw [line width = 1] (1,4)--(2,4);
\draw [line width = 1] (1,5)--(2,5);

\draw [line width = 1] (3,4)--(3,5);
\draw [line width = 1] (4,4)--(4,5);
\draw [line width = 1] (5,4)--(5,5);
\draw [line width = 1] (6,4)--(6,5);

\node at (1.5,1.5) {$+$};
\node at (1.5,2.5) {$+$};
\node at (1.5,3.5) {$+$};
\node at (1.5,4.5) {$0$};
\node at (1.5,5.5) {$0$};

\node at (2.5,4.5) {$-$};
\node at (3.5,4.5) {$-$};
\node at (4.5,4.5) {$-$};
\node at (5.5,4.5) {$-$};
\node at (6.5,4.5) {$-$};

\node at (.5,.5) {$-$};
\node at (.5,1.5) {$-$};
\node at (.5,2.5) {$-$};
\node at (.5,3.5) {$-$};
\node at (.5,4.5) {$-$};
\node at (.5,5.5) {$0$};

\node at (1.5,.5) {$0$};
\end{tikzpicture}
\caption{Step 3 (-)}
\end{subfigure}
~~~~
\begin{subfigure}{.16\textwidth}
\begin{tikzpicture}[scale=.33]

\draw [line width = 1,<->] (0,7)--(0,0)--(8,0);
\draw [line width = 2] (1,5)--(2,5)--(2,4)--(7,4)--(7,3)--(1,3)--(1,5);
\draw [line width = 1] (1,6)--(2,6)--(2,5)--(7,5)--(7,4)--(1,4)--(1,6);
\draw [line width = 1] (1,0)--(2,0)--(2,6)--(0,6)--(0,5)--(1,5)--(1,0);
\draw [line width = 1] (0,0)--(2,0)--(2,1)--(1,1)--(1,6)--(0,6)--(0,0);
\draw [line width = 1] (1,0)--(1,1);
\draw [line width = 1] (0,1)--(1,1);
\draw [line width = 1] (0,2)--(1,2);
\draw [line width = 1] (0,3)--(1,3);
\draw [line width = 1] (0,4)--(1,4);
\draw [line width = 1] (0,5)--(1,5);

\draw [line width = 1] (1,2)--(2,2);
\draw [line width = 1] (1,3)--(2,3);
\draw [line width = 1] (1,4)--(2,4);
\draw [line width = 1] (1,5)--(2,5);

\draw [line width = 1] (3,4)--(3,5);
\draw [line width = 1] (4,4)--(4,5);
\draw [line width = 1] (5,4)--(5,5);
\draw [line width = 1] (6,4)--(6,5);

\draw [line width = 1] (3,3)--(3,4);
\draw [line width = 1] (4,3)--(4,4);
\draw [line width = 1] (5,3)--(5,4);
\draw [line width = 1] (6,3)--(6,4);

\node at (1.5,1.5) {$+$};
\node at (1.5,2.5) {$+$};
\node at (1.5,3.5) {$0$};
\node at (1.5,4.5) {$-$};
\node at (1.5,5.5) {$0$};

\node at (2.5,4.5) {$-$};
\node at (3.5,4.5) {$-$};
\node at (4.5,4.5) {$-$};
\node at (5.5,4.5) {$-$};
\node at (6.5,4.5) {$-$};

\node at (2.5,3.5) {$-$};
\node at (3.5,3.5) {$-$};
\node at (4.5,3.5) {$-$};
\node at (5.5,3.5) {$-$};
\node at (6.5,3.5) {$-$};

\node at (.5,.5) {$-$};
\node at (.5,1.5) {$-$};
\node at (.5,2.5) {$-$};
\node at (.5,3.5) {$-$};
\node at (.5,4.5) {$-$};
\node at (.5,5.5) {$0$};

\node at (1.5,.5) {$0$};
\end{tikzpicture}
\caption{Step 4 (-)}
\end{subfigure}
~~~~
\begin{subfigure}{.16\textwidth}
\begin{tikzpicture}[scale=.33]

\draw [line width = 1,<->] (0,7)--(0,0)--(8,0);
\draw [line width = 2] (1,4)--(2,4)--(2,3)--(7,3)--(7,2)--(1,2)--(1,4);
\draw [line width = 1] (1,5)--(2,5)--(2,4)--(7,4)--(7,3)--(1,3)--(1,5);
\draw [line width = 1] (1,6)--(2,6)--(2,5)--(7,5)--(7,4)--(1,4)--(1,6);
\draw [line width = 1] (1,0)--(2,0)--(2,6)--(0,6)--(0,5)--(1,5)--(1,0);
\draw [line width = 1] (0,0)--(2,0)--(2,1)--(1,1)--(1,6)--(0,6)--(0,0);
\draw [line width = 1] (1,0)--(1,1);
\draw [line width = 1] (0,1)--(1,1);
\draw [line width = 1] (0,2)--(1,2);
\draw [line width = 1] (0,3)--(1,3);
\draw [line width = 1] (0,4)--(1,4);
\draw [line width = 1] (0,5)--(1,5);

\draw [line width = 1] (1,2)--(2,2);
\draw [line width = 1] (1,3)--(2,3);
\draw [line width = 1] (1,4)--(2,4);
\draw [line width = 1] (1,5)--(2,5);

\draw [line width = 1] (3,4)--(3,5);
\draw [line width = 1] (4,4)--(4,5);
\draw [line width = 1] (5,4)--(5,5);
\draw [line width = 1] (6,4)--(6,5);

\draw [line width = 1] (3,3)--(3,4);
\draw [line width = 1] (4,3)--(4,4);
\draw [line width = 1] (5,3)--(5,4);
\draw [line width = 1] (6,3)--(6,4);

\draw [line width = 1] (3,3)--(3,2);
\draw [line width = 1] (4,3)--(4,2);
\draw [line width = 1] (5,3)--(5,2);
\draw [line width = 1] (6,3)--(6,2);

\node at (1.5,1.5) {$+$};
\node at (1.5,2.5) {$0$};
\node at (1.5,3.5) {$-$};
\node at (1.5,4.5) {$-$};
\node at (1.5,5.5) {$0$};

\node at (2.5,4.5) {$-$};
\node at (3.5,4.5) {$-$};
\node at (4.5,4.5) {$-$};
\node at (5.5,4.5) {$-$};
\node at (6.5,4.5) {$-$};

\node at (2.5,3.5) {$-$};
\node at (3.5,3.5) {$-$};
\node at (4.5,3.5) {$-$};
\node at (5.5,3.5) {$-$};
\node at (6.5,3.5) {$-$};

\node at (2.5,2.5) {$-$};
\node at (3.5,2.5) {$-$};
\node at (4.5,2.5) {$-$};
\node at (5.5,2.5) {$-$};
\node at (6.5,2.5) {$-$};

\node at (.5,.5) {$-$};
\node at (.5,1.5) {$-$};
\node at (.5,2.5) {$-$};
\node at (.5,3.5) {$-$};
\node at (.5,4.5) {$-$};
\node at (.5,5.5) {$0$};

\node at (1.5,.5) {$0$};
\end{tikzpicture}
\caption{Step 5 (-)}
\end{subfigure}
~~~~
\begin{subfigure}{.16\textwidth}
\begin{tikzpicture}[scale=.33]

\draw [line width = 1,<->] (0,7)--(0,0)--(8,0);
\draw [line width = 2] (1,3)--(2,3)--(2,2)--(7,2)--(7,1)--(1,1)--(1,3);
\draw [line width = 1] (1,4)--(2,4)--(2,3)--(7,3)--(7,2)--(1,2)--(1,4);
\draw [line width = 1] (1,5)--(2,5)--(2,4)--(7,4)--(7,3)--(1,3)--(1,5);
\draw [line width = 1] (1,6)--(2,6)--(2,5)--(7,5)--(7,4)--(1,4)--(1,6);
\draw [line width = 1] (1,0)--(2,0)--(2,6)--(0,6)--(0,5)--(1,5)--(1,0);
\draw [line width = 1] (0,0)--(2,0)--(2,1)--(1,1)--(1,6)--(0,6)--(0,0);
\draw [line width = 1] (1,0)--(1,1);
\draw [line width = 1] (0,1)--(1,1);
\draw [line width = 1] (0,2)--(1,2);
\draw [line width = 1] (0,3)--(1,3);
\draw [line width = 1] (0,4)--(1,4);
\draw [line width = 1] (0,5)--(1,5);

\draw [line width = 1] (1,2)--(2,2);
\draw [line width = 1] (1,3)--(2,3);
\draw [line width = 1] (1,4)--(2,4);
\draw [line width = 1] (1,5)--(2,5);

\draw [line width = 1] (3,4)--(3,5);
\draw [line width = 1] (4,4)--(4,5);
\draw [line width = 1] (5,4)--(5,5);
\draw [line width = 1] (6,4)--(6,5);

\draw [line width = 1] (3,3)--(3,4);
\draw [line width = 1] (4,3)--(4,4);
\draw [line width = 1] (5,3)--(5,4);
\draw [line width = 1] (6,3)--(6,4);

\draw [line width = 1] (3,3)--(3,2);
\draw [line width = 1] (4,3)--(4,2);
\draw [line width = 1] (5,3)--(5,2);
\draw [line width = 1] (6,3)--(6,2);

\draw [line width = 1] (3,1)--(3,2);
\draw [line width = 1] (4,1)--(4,2);
\draw [line width = 1] (5,1)--(5,2);
\draw [line width = 1] (6,1)--(6,2);

\node at (1.5,1.5) {$0$};
\node at (1.5,2.5) {$-$};
\node at (1.5,3.5) {$-$};
\node at (1.5,4.5) {$-$};
\node at (1.5,5.5) {$0$};

\node at (2.5,4.5) {$-$};
\node at (3.5,4.5) {$-$};
\node at (4.5,4.5) {$-$};
\node at (5.5,4.5) {$-$};
\node at (6.5,4.5) {$-$};

\node at (2.5,3.5) {$-$};
\node at (3.5,3.5) {$-$};
\node at (4.5,3.5) {$-$};
\node at (5.5,3.5) {$-$};
\node at (6.5,3.5) {$-$};

\node at (2.5,2.5) {$-$};
\node at (3.5,2.5) {$-$};
\node at (4.5,2.5) {$-$};
\node at (5.5,2.5) {$-$};
\node at (6.5,2.5) {$-$};

\node at (2.5,1.5) {$-$};
\node at (3.5,1.5) {$-$};
\node at (4.5,1.5) {$-$};
\node at (5.5,1.5) {$-$};
\node at (6.5,1.5) {$-$};

\node at (.5,.5) {$-$};
\node at (.5,1.5) {$-$};
\node at (.5,2.5) {$-$};
\node at (.5,3.5) {$-$};
\node at (.5,4.5) {$-$};
\node at (.5,5.5) {$0$};

\node at (1.5,.5) {$0$};
\end{tikzpicture}
\caption{Step 6 (-)}
\end{subfigure}
~~~~
\begin{subfigure}{.16\textwidth}
\begin{tikzpicture}[scale=.33]

\draw [line width = 1,<->] (0,7)--(0,0)--(8,0);
\draw [line width = 2] (0,5)--(1,5)--(1,4)--(6,4)--(6,3)--(0,3)--(0,5);
\draw [line width = 1] (1,3)--(2,3)--(2,2)--(7,2)--(7,1)--(1,1)--(1,3);
\draw [line width = 1] (1,4)--(2,4)--(2,3)--(7,3)--(7,2)--(1,2)--(1,4);
\draw [line width = 1] (1,5)--(2,5)--(2,4)--(7,4)--(7,3)--(1,3)--(1,5);
\draw [line width = 1] (1,6)--(2,6)--(2,5)--(7,5)--(7,4)--(1,4)--(1,6);
\draw [line width = 1] (1,0)--(2,0)--(2,6)--(0,6)--(0,5)--(1,5)--(1,0);
\draw [line width = 1] (0,0)--(2,0)--(2,1)--(1,1)--(1,6)--(0,6)--(0,0);
\draw [line width = 1] (1,0)--(1,1);
\draw [line width = 1] (0,1)--(1,1);
\draw [line width = 1] (0,2)--(1,2);
\draw [line width = 1] (0,3)--(1,3);
\draw [line width = 1] (0,4)--(1,4);
\draw [line width = 1] (0,5)--(1,5);

\draw [line width = 1] (1,2)--(2,2);
\draw [line width = 1] (1,3)--(2,3);
\draw [line width = 1] (1,4)--(2,4);
\draw [line width = 1] (1,5)--(2,5);

\draw [line width = 1] (3,4)--(3,5);
\draw [line width = 1] (4,4)--(4,5);
\draw [line width = 1] (5,4)--(5,5);
\draw [line width = 1] (6,4)--(6,5);

\draw [line width = 1] (3,3)--(3,4);
\draw [line width = 1] (4,3)--(4,4);
\draw [line width = 1] (5,3)--(5,4);
\draw [line width = 1] (6,3)--(6,4);

\draw [line width = 1] (3,3)--(3,2);
\draw [line width = 1] (4,3)--(4,2);
\draw [line width = 1] (5,3)--(5,2);
\draw [line width = 1] (6,3)--(6,2);

\draw [line width = 1] (3,1)--(3,2);
\draw [line width = 1] (4,1)--(4,2);
\draw [line width = 1] (5,1)--(5,2);
\draw [line width = 1] (6,1)--(6,2);

\node at (1.5,1.5) {$0$};
\node at (1.5,2.5) {$-$};
\node at (1.5,3.5) {$0$};
\node at (1.5,4.5) {$-$};
\node at (1.5,5.5) {$0$};

\node at (2.5,4.5) {$-$};
\node at (3.5,4.5) {$-$};
\node at (4.5,4.5) {$-$};
\node at (5.5,4.5) {$-$};
\node at (6.5,4.5) {$-$};

\node at (2.5,3.5) {$0$};
\node at (3.5,3.5) {$0$};
\node at (4.5,3.5) {$0$};
\node at (5.5,3.5) {$0$};
\node at (6.5,3.5) {$-$};

\node at (2.5,2.5) {$-$};
\node at (3.5,2.5) {$-$};
\node at (4.5,2.5) {$-$};
\node at (5.5,2.5) {$-$};
\node at (6.5,2.5) {$-$};

\node at (2.5,1.5) {$-$};
\node at (3.5,1.5) {$-$};
\node at (4.5,1.5) {$-$};
\node at (5.5,1.5) {$-$};
\node at (6.5,1.5) {$-$};

\node at (.5,.5) {$-$};
\node at (.5,1.5) {$-$};
\node at (.5,2.5) {$-$};
\node at (.5,3.5) {$0$};
\node at (.5,4.5) {$0$};
\node at (.5,5.5) {$0$};

\node at (1.5,.5) {$0$};
\end{tikzpicture}
\caption{Step 7 (+)}
\end{subfigure}
~~~~
\begin{subfigure}{.16\textwidth}
\begin{tikzpicture}[scale=.33]

\draw [line width = 1,<->] (0,7)--(0,0)--(8,0);
\draw [line width = 2] (1,4)--(1,5)--(7,5)--(7,3)--(6,3)--(6,4)--(1,4);
\draw [line width = 1] (0,5)--(1,5)--(1,4)--(6,4)--(6,3)--(0,3)--(0,5);
\draw [line width = 1] (1,3)--(2,3)--(2,2)--(7,2)--(7,1)--(1,1)--(1,3);
\draw [line width = 1] (1,4)--(2,4)--(2,3)--(7,3)--(7,2)--(1,2)--(1,4);
\draw [line width = 1] (1,5)--(2,5)--(2,4)--(7,4)--(7,3)--(1,3)--(1,5);
\draw [line width = 1] (1,6)--(2,6)--(2,5)--(7,5)--(7,4)--(1,4)--(1,6);
\draw [line width = 1] (1,0)--(2,0)--(2,6)--(0,6)--(0,5)--(1,5)--(1,0);
\draw [line width = 1] (0,0)--(2,0)--(2,1)--(1,1)--(1,6)--(0,6)--(0,0);
\draw [line width = 1] (1,0)--(1,1);
\draw [line width = 1] (0,1)--(1,1);
\draw [line width = 1] (0,2)--(1,2);
\draw [line width = 1] (0,3)--(1,3);
\draw [line width = 1] (0,4)--(1,4);
\draw [line width = 1] (0,5)--(1,5);

\draw [line width = 1] (1,2)--(2,2);
\draw [line width = 1] (1,3)--(2,3);
\draw [line width = 1] (1,4)--(2,4);
\draw [line width = 1] (1,5)--(2,5);

\draw [line width = 1] (3,4)--(3,5);
\draw [line width = 1] (4,4)--(4,5);
\draw [line width = 1] (5,4)--(5,5);
\draw [line width = 1] (6,4)--(6,5);

\draw [line width = 1] (3,3)--(3,4);
\draw [line width = 1] (4,3)--(4,4);
\draw [line width = 1] (5,3)--(5,4);
\draw [line width = 1] (6,3)--(6,4);

\draw [line width = 1] (3,3)--(3,2);
\draw [line width = 1] (4,3)--(4,2);
\draw [line width = 1] (5,3)--(5,2);
\draw [line width = 1] (6,3)--(6,2);

\draw [line width = 1] (3,1)--(3,2);
\draw [line width = 1] (4,1)--(4,2);
\draw [line width = 1] (5,1)--(5,2);
\draw [line width = 1] (6,1)--(6,2);

\node at (1.5,1.5) {$0$};
\node at (1.5,2.5) {$-$};
\node at (1.5,3.5) {$0$};
\node at (1.5,4.5) {$0$};
\node at (1.5,5.5) {$0$};

\node at (2.5,4.5) {$0$};
\node at (3.5,4.5) {$0$};
\node at (4.5,4.5) {$0$};
\node at (5.5,4.5) {$0$};
\node at (6.5,4.5) {$0$};

\node at (2.5,3.5) {$0$};
\node at (3.5,3.5) {$0$};
\node at (4.5,3.5) {$0$};
\node at (5.5,3.5) {$0$};
\node at (6.5,3.5) {$0$};

\node at (2.5,2.5) {$-$};
\node at (3.5,2.5) {$-$};
\node at (4.5,2.5) {$-$};
\node at (5.5,2.5) {$-$};
\node at (6.5,2.5) {$-$};

\node at (2.5,1.5) {$-$};
\node at (3.5,1.5) {$-$};
\node at (4.5,1.5) {$-$};
\node at (5.5,1.5) {$-$};
\node at (6.5,1.5) {$-$};

\node at (.5,.5) {$-$};
\node at (.5,1.5) {$-$};
\node at (.5,2.5) {$-$};
\node at (.5,3.5) {$0$};
\node at (.5,4.5) {$0$};
\node at (.5,5.5) {$0$};

\node at (1.5,.5) {$0$};
\end{tikzpicture}
\caption{Step 8 (+)}
\end{subfigure}
~~~~
\begin{subfigure}{.16\textwidth}
\begin{tikzpicture}[scale=.33]

\draw [line width = 1,<->] (0,7)--(0,0)--(8,0);
\draw [line width = 2] (0,3)--(1,3)--(1,2)--(6,2)--(6,1)--(0,1)--(0,3);
\draw [line width = 1] (1,4)--(1,5)--(7,5)--(7,3)--(6,3)--(6,4)--(1,4);
\draw [line width = 1] (0,5)--(1,5)--(1,4)--(6,4)--(6,3)--(0,3)--(0,5);
\draw [line width = 1] (1,3)--(2,3)--(2,2)--(7,2)--(7,1)--(1,1)--(1,3);
\draw [line width = 1] (1,4)--(2,4)--(2,3)--(7,3)--(7,2)--(1,2)--(1,4);
\draw [line width = 1] (1,5)--(2,5)--(2,4)--(7,4)--(7,3)--(1,3)--(1,5);
\draw [line width = 1] (1,6)--(2,6)--(2,5)--(7,5)--(7,4)--(1,4)--(1,6);
\draw [line width = 1] (1,0)--(2,0)--(2,6)--(0,6)--(0,5)--(1,5)--(1,0);
\draw [line width = 1] (0,0)--(2,0)--(2,1)--(1,1)--(1,6)--(0,6)--(0,0);
\draw [line width = 1] (1,0)--(1,1);
\draw [line width = 1] (0,1)--(1,1);
\draw [line width = 1] (0,2)--(1,2);
\draw [line width = 1] (0,3)--(1,3);
\draw [line width = 1] (0,4)--(1,4);
\draw [line width = 1] (0,5)--(1,5);

\draw [line width = 1] (1,2)--(2,2);
\draw [line width = 1] (1,3)--(2,3);
\draw [line width = 1] (1,4)--(2,4);
\draw [line width = 1] (1,5)--(2,5);

\draw [line width = 1] (3,4)--(3,5);
\draw [line width = 1] (4,4)--(4,5);
\draw [line width = 1] (5,4)--(5,5);
\draw [line width = 1] (6,4)--(6,5);

\draw [line width = 1] (3,3)--(3,4);
\draw [line width = 1] (4,3)--(4,4);
\draw [line width = 1] (5,3)--(5,4);
\draw [line width = 1] (6,3)--(6,4);

\draw [line width = 1] (3,3)--(3,2);
\draw [line width = 1] (4,3)--(4,2);
\draw [line width = 1] (5,3)--(5,2);
\draw [line width = 1] (6,3)--(6,2);

\draw [line width = 1] (3,1)--(3,2);
\draw [line width = 1] (4,1)--(4,2);
\draw [line width = 1] (5,1)--(5,2);
\draw [line width = 1] (6,1)--(6,2);

\node at (1.5,1.5) {$+$};
\node at (1.5,2.5) {$-$};
\node at (1.5,3.5) {$0$};
\node at (1.5,4.5) {$0$};
\node at (1.5,5.5) {$0$};

\node at (2.5,4.5) {$0$};
\node at (3.5,4.5) {$0$};
\node at (4.5,4.5) {$0$};
\node at (5.5,4.5) {$0$};
\node at (6.5,4.5) {$0$};

\node at (2.5,3.5) {$0$};
\node at (3.5,3.5) {$0$};
\node at (4.5,3.5) {$0$};
\node at (5.5,3.5) {$0$};
\node at (6.5,3.5) {$0$};

\node at (2.5,2.5) {$-$};
\node at (3.5,2.5) {$-$};
\node at (4.5,2.5) {$-$};
\node at (5.5,2.5) {$-$};
\node at (6.5,2.5) {$-$};

\node at (2.5,1.5) {$0$};
\node at (3.5,1.5) {$0$};
\node at (4.5,1.5) {$0$};
\node at (5.5,1.5) {$0$};
\node at (6.5,1.5) {$-$};

\node at (.5,.5) {$-$};
\node at (.5,1.5) {$0$};
\node at (.5,2.5) {$0$};
\node at (.5,3.5) {$0$};
\node at (.5,4.5) {$0$};
\node at (.5,5.5) {$0$};

\node at (1.5,.5) {$0$};
\end{tikzpicture}
\caption{Step 9 (+)}
\end{subfigure}
~~
\begin{subfigure}{.16\textwidth}
\begin{tikzpicture}[scale=.33]

\draw [line width = 1,<->] (0,7)--(0,0)--(8,0);
\draw [line width = 1, fill=lightgray] (0,0)--(1,0)--(1,1)--(0,1)--(0,0);
\draw [line width = 1, fill=lightgray] (1,1)--(2,1)--(2,2)--(1,2)--(1,1);
\draw [line width = 2] (1,2)--(1,3)--(7,3)--(7,1)--(6,1)--(6,2)--(1,2);
\draw [line width = 1] (0,3)--(1,3)--(1,2)--(6,2)--(6,1)--(0,1)--(0,3);
\draw [line width = 1] (1,4)--(1,5)--(7,5)--(7,3)--(6,3)--(6,4)--(1,4);
\draw [line width = 1] (0,5)--(1,5)--(1,4)--(6,4)--(6,3)--(0,3)--(0,5);
\draw [line width = 1] (1,3)--(2,3)--(2,2)--(7,2)--(7,1)--(1,1)--(1,3);
\draw [line width = 1] (1,4)--(2,4)--(2,3)--(7,3)--(7,2)--(1,2)--(1,4);
\draw [line width = 1] (1,5)--(2,5)--(2,4)--(7,4)--(7,3)--(1,3)--(1,5);
\draw [line width = 1] (1,6)--(2,6)--(2,5)--(7,5)--(7,4)--(1,4)--(1,6);
\draw [line width = 1] (1,0)--(2,0)--(2,6)--(0,6)--(0,5)--(1,5)--(1,0);
\draw [line width = 1] (0,0)--(2,0)--(2,1)--(1,1)--(1,6)--(0,6)--(0,0);
\draw [line width = 1] (1,0)--(1,1);
\draw [line width = 1] (0,1)--(1,1);
\draw [line width = 1] (0,2)--(1,2);
\draw [line width = 1] (0,3)--(1,3);
\draw [line width = 1] (0,4)--(1,4);
\draw [line width = 1] (0,5)--(1,5);

\draw [line width = 1] (1,2)--(2,2);
\draw [line width = 1] (1,3)--(2,3);
\draw [line width = 1] (1,4)--(2,4);
\draw [line width = 1] (1,5)--(2,5);

\draw [line width = 1] (3,4)--(3,5);
\draw [line width = 1] (4,4)--(4,5);
\draw [line width = 1] (5,4)--(5,5);
\draw [line width = 1] (6,4)--(6,5);

\draw [line width = 1] (3,3)--(3,4);
\draw [line width = 1] (4,3)--(4,4);
\draw [line width = 1] (5,3)--(5,4);
\draw [line width = 1] (6,3)--(6,4);

\draw [line width = 1] (3,3)--(3,2);
\draw [line width = 1] (4,3)--(4,2);
\draw [line width = 1] (5,3)--(5,2);
\draw [line width = 1] (6,3)--(6,2);

\draw [line width = 1] (3,1)--(3,2);
\draw [line width = 1] (4,1)--(4,2);
\draw [line width = 1] (5,1)--(5,2);
\draw [line width = 1] (6,1)--(6,2);

\node at (1.5,1.5) {$+$};
\node at (1.5,2.5) {$0$};
\node at (1.5,3.5) {$0$};
\node at (1.5,4.5) {$0$};
\node at (1.5,5.5) {$0$};

\node at (2.5,4.5) {$0$};
\node at (3.5,4.5) {$0$};
\node at (4.5,4.5) {$0$};
\node at (5.5,4.5) {$0$};
\node at (6.5,4.5) {$0$};

\node at (2.5,3.5) {$0$};
\node at (3.5,3.5) {$0$};
\node at (4.5,3.5) {$0$};
\node at (5.5,3.5) {$0$};
\node at (6.5,3.5) {$0$};

\node at (2.5,2.5) {$0$};
\node at (3.5,2.5) {$0$};
\node at (4.5,2.5) {$0$};
\node at (5.5,2.5) {$0$};
\node at (6.5,2.5) {$0$};

\node at (2.5,1.5) {$0$};
\node at (3.5,1.5) {$0$};
\node at (4.5,1.5) {$0$};
\node at (5.5,1.5) {$0$};
\node at (6.5,1.5) {$0$};

\node at (.5,.5) {$-$};
\node at (.5,1.5) {$0$};
\node at (.5,2.5) {$0$};
\node at (.5,3.5) {$0$};
\node at (.5,4.5) {$0$};
\node at (.5,5.5) {$0$};

\node at (1.5,.5) {$0$};
\end{tikzpicture}
\caption{Step 10 (+)}
\end{subfigure}
~~~~
\caption{The polynomial $B_3(7)$ is generated by $\{G_1(7),G_2(7),G_3(7),G_4(7)\}$.}
\label{fig:special-form-33}
\end{figure}

\end{pf}

\begin{propo} $\{B_1(k), B_2(k), B_3(k)\}$ and $\{G_1(k), G_2(k), G_3(k), G_4(k)\}$ generate the same ideal in $\mathbb{Z}[X,Y]$.
\end{propo}

\begin{pf} This follows from Propositions~\ref{p:prop5},~\ref{p:prop6}.
\end{pf}

\begin{propo} One has the following $D$-reductions
\begin{equation}
\begin{aligned}
S(B_1(k),B_2(k))&=-y^kB_1(k)+x^{k-1}B_2(k)+\left (x^{k-1}\cdot \frac{y^k-1}{y-1}-y^k\cdot\frac{x^{k-1}-1}{x-1}\right )\cdot B_3(k)\\
S(B_1(k),B_3(k))&=B_2(k)+\frac{y^k-1}{y-1}\cdot B_3(k)\\
S(B_2(k),B_3(k))&=B_1(k)+\frac{x^{k-1}-1}{x-1}\cdot B_3(k).
\end{aligned}
\end{equation}
Consequently, $\{B_1(k), B_2(k), B_3(k)\}$ is a Gr\"obner basis.
\end{propo}

\begin{pf} The leading monomial of $B_1(k)$ is $y^{k+1}$, the leading monomial of $B_2(k)$ is $x^k$ and the leading monomial of $B_3(k)$ is $xy$.
We reduce the $S$-polynomials related to the set $\{B_1(k),B_2(k),B_3(k)\}$:
\begin{equation}
\begin{aligned}
S(B_1(k),B_2(k))&=x^k\cdot B_1(k)-y^{k+1}\cdot B_2(k)\\
&=x^k\cdot \left ( \frac{y^{k+2}-1}{y-1}+x\cdot\frac{x^{k-1}-1}{x-1}\right )-y^{k+1}\cdot \left ( \frac{y^{k+1}-1}{y-1}+x\cdot\frac{x^{k}-1}{x-1}\right )\\
&=\frac{-xy^{2k+2}-x^ky^{k+2}+x^k+y^{2k+2}-y^{k+1}+x^{2k}y-x^{k+1}y+xy^{k+2}-x^{2k}+y^{k+1}x^{k+1}}{(x-1)(y-1)}\\
&=-y^kB_1(k)+x^{k-1}B_2(k)+\left (x^{k-1}\cdot \frac{y^k-1}{y-1}-y^k\cdot\frac{x^{k-1}-1}{x-1}\right )\cdot B_3(k).
\end{aligned}
\end{equation}

\begin{equation}
\begin{aligned}
S(B_1(k),B_3(k))&=x\cdot B_1(k)-y^{k}\cdot B_3(k)\\
&=x\cdot \left ( \frac{y^{k+2}-1}{y-1}+x\cdot\frac{x^{k-1}-1}{x-1}\right )-y^{k}\cdot (xy-1)\\
&=\frac{x+x^{k+1}y-x^2y-x^{k+1}+x^2y^{k+1}-xy^{k+1}+xy^{k+1}-xy^k-y^{k+1}+y^k}{(x-1)(y-1)}\\
&=B_2(k)+\frac{y^k-1}{y-1}\cdot B_3(k).
\end{aligned}
\end{equation}

\begin{equation}
\begin{aligned}
S(B_2(k),B_3(k))&=y\cdot B_2(k)-x^{k-1}\cdot B_3(k)\\
&=y\cdot \left ( \frac{y^{k+1}-1}{y-1}+x\cdot\frac{x^{k}-1}{x-1}\right )-x^{k-1}\cdot (xy-1)\\
&=\frac{xy^{k+2}-y^{k+2}-x^k+x^ky^2-yx^{k-1}-xy^2+y+x^{k-1}}{(x-1)(y-1)}\\
&=B_1(k)+\frac{x^{k-1}-1}{x-1}\cdot B_3(k).
\end{aligned}
\end{equation}

We show now that all above reductions are $D$-reductions by looking at the elimination of the terms of highest degree in the $S$-polynomials.

The terms of highest degrees in $S(B_1(k),B_2(k))$, after the initial reduction (underlined below)
\begin{equation}
\begin{gathered}
x^k\cdot B_1(k)-y^{k+1}\cdot B_2(k)=x^k(\underline{{\bf y^{k+1}}}+y^k+y^{k-1}+\cdots+x^{k-1}+x^{k-2}+x^{k-3}+\cdots)\\
-y^{k+1}(y^k+y^{k-1}+y^{k-2}+\cdots+\underline{{\bf x^k}}+x^{k-1}+x^{k-2}+\cdots),
\end{gathered}
\end{equation} are (in this order) $$-y^{2k+1}+x^ky^k-x^{k-1}y^{k+1}-y^{2k}.$$

The terms $-y^{2k+1}-y^{2k}$ are contained in
\begin{equation}
-y^kB_1(k)=-y^k(\underline{\bf y^{k+1}+y^k}+y^{k-1}+\cdots +x^{k-1}+x^{k-2}+x^{k-3}+\cdots),
\end{equation} which does not contains terms of higher degree then $x^ky^k-x^{k-1}y^{k+1}.$

The remaining terms $x^ky^k-x^{k-1}y^{k+1}$ are contained in
\begin{equation}
\begin{gathered}
\left (x^{k-1}\cdot \frac{y^k-1}{y-1}-y^k\cdot\frac{x^{k-1}-1}{x-1}\right )\cdot B_3(k)\\
=\left [x^{k-1}(\underline{{\bf y^{k-1}}}+y^{k-2}+y^{k-3}+\cdots)-y^k(\underline{{\bf x^{k-2}}}+x^{k-3}+\cdots) \right ](xy-1),
\end{gathered}
\end{equation}
which also does not contain terms of higher degree then $x^ky^k-x^{k-1}y^{k+1}.$

The term of highest degrees in $S(B_1(k),B_3(k))$, after the initial reduction (underlined below)
\begin{equation}
\begin{gathered}
x\cdot B_1(k)-y^{k}\cdot B_3(k)=x(\underline{{\bf y^{k+1}}}+y^k+y^{k-1}+\cdots+x^{k-1}+x^{k-2}+x^{k-3}+\cdots)-y^k(\underline{{\bf xy}}-1)
\end{gathered}
\end{equation}
is $xy^k$. This term is contained in
\begin{equation}
\frac{y^k-1}{y-1}\cdot B_3(k)=(y^{k-1}+y^{k-2}+\cdots)(xy-1),
\end{equation} which does not contain terms of higher degree then $xy^k$.

The term of highest degrees in $S(B_2(k),B_3(k))$, after the initial reduction (underlined below)
\begin{equation}
\begin{gathered}
y\cdot B_2(k)-x^{k-1}\cdot B_3(k)=y(y^k+y^{k-1}+y^{k-2}+\cdots+\underline{{\bf x^k}}+x^{k-1}+x^{k-2}+\cdots)-x^{k-1}(\underline{{\bf xy}}-1)
\end{gathered}
\end{equation} is $y^{k+1}$. This term is contained in $B_1(k)$, which does not contain terms of higher degree then $y^{k+1}$.

As all higher coefficients are equal to $1$, we do not need to consider the $G$-polynomials.
\end{pf}

\section{Proof of Theorem~\ref{thm-main}}

Consider a $q\times p, q\ge p\ge 1,$ rectangle. Using the presence of $B_3(k)$ in the Gr\" obner basis, and Theorem~\ref{thm-groebner-tiling}, the existence of a signed tiling becomes equivalent to deciding when the polynomial:
\begin{equation}
P_{p,q}(x)=1+2x+3x^2+\dots +px^{p-1}+px^{p}+\dots +px^{q-1}+(p-1)x^q+(p-2)x^{q+1} +\dots +2x^{p+q-3}+x^{p+q-2}
\end{equation}
is divisible by the polynomial:
\begin{equation}
Q(x)=1+x+x^2+\dots +x^{n-1}.
\end{equation}

If $p+q-1<n$, then $\deg Q>\deg P_{p,q}$, so divisibility does not hold. If $p+q-1\ge n$, we look at $P_{p,q}$ as a sum of $p$ polynomials with all coefficients equal to 1:
\begin{equation}
\begin{aligned}
P_{p,q}(x)=1+x+x^2+x^3+\dots +x^{p-1}+x^{p}+\dots &+x^{q-1}+x^q+x^{q+1} +\dots +x^{p+q-4}+x^{p+q-3}+x^{p+q-2}\\
          +x+x^2+x^3+\dots +x^{p-1}+x^{p}+\dots &+x^{q-1}+x^q+x^{q+1} +x^{p+q-4}+\dots +x^{p+q-3}\\
          +x^3+\dots +x^{p-1}+x^{p}+\dots &+x^{q-1}+x^q+x^{q+1} +\dots +x^{p+q-4}\\
          \dots \dots &\dots \dots\\
        +x^{p}+\dots &+x^{q-1}.
\end{aligned}
\end{equation}

Assume that $p+q-1=nm+r, 0\le r<n,$ and $p=ns+t, 0\le t<n.$ The remainder $R_{p,q}(x)$ of the division of $P_{p,q}(x)$ by $Q(x)$ is the sum of the remainders of the division of the $p$ polynomials above by $Q(x)$.

If $r$ is odd, one has the following sequence of remainders, each remainder written in a separate pair of parentheses:
\begin{equation}
\begin{aligned}
R_{p,q}(x)=&(1+x+x^2+\dots +x^{r-1})\\
+&(x+x^2+\dots +x^{r-2})\\
+&(x^2+\dots +x^{r-3})\\
&\ \ \ \dots \dots \dots\\
+&(x^{\frac{r-1}{2}}+x^{\frac{r+1}{2}})+(0)-(x^{\frac{r-1}{2}}+x^{\frac{r+1}{2}})\\
&\ \ \ \dots \dots \dots\\
-&(x+x^2+\dots +x^{r-2})\\
-&(1+x+x^2+\dots +x^{r-1})\\
+&(x^{r+1}+x^{r+3}+\dots +x^{n-3}+x^{n-2})\\
+&(x^{r+2}+\dots +x^{n-3})\\
&\ \ \ \dots \dots \dots\\
+&(x^{\frac{r+n}{2}})-(x^{\frac{r+n}{2}})\\
&\ \ \ \dots \dots \dots\\
-&(x^{r+2}+\dots +x^{n-3})\\
-&(x^{r+1}+x^{r+3}+\dots +x^{n-3}+x^{n-2})\\
&\ \ \ \dots \dots \dots\\
\end{aligned}
\end{equation}

If $p\ge n$, the sequence of remainders above is periodic with period $n$, given by the part of the sequence shown above, and the sum of any subsequence of $n$ consecutive remainders is 0. So if $p$ is divisible by $n$, $P_{p,q}(x)$ is divisible by $Q(x)$. If $p$ is not divisible by $n$, then doing first the cancellation as above and then using the symmetry of the sequence of remainders about the remainder equal to 0, the sum of the sequence of remainders equals 0 only if $r+1=t$, that is, only if $q$ is divisible by $n$.

If $r$ is even, one has the following sequence of remainders, each remainder written in a separate pair of parentheses:
\begin{equation}
\begin{aligned}
R_{p,q}(x)=&(1+x+x^2+\dots +x^{r-1})\\
+&(x+x^2+\dots +x^{r-2})\\
+&(x^2+\dots +x^{r-3})\\
&\ \ \ \dots \dots \dots\\
+&(x^{\frac{r}{2}})-(x^{\frac{r}{2}})\\
&\ \ \ \dots \dots \dots\\
-&((x+x^2+\dots +x^{r-2}))\\
-&(1+x+x^2+\dots +x^{r-1})\\
+&(x^{r+1}+x^{r+2}+\dots +x^{n-3}+x^{n-2})\\
+&(x^{r+2}+\dots +x^{n-3})\\
&\ \ \ \dots \dots \dots\\
+&(x^{\frac{r+n-1}{2}}+x^{\frac{r+n+1}{2}})+(0)-(x^{\frac{r+n-1}{2}}+x^{\frac{r+n+1}{2}})\\
&\ \ \ \dots \dots \dots\\
-&(x^{r+2}+\dots +x^{n-3})\\
-&(x^{r+1}+x^{r+2}+\dots +x^{n-3}+x^{n-2})\\
&\ \ \ \dots \dots \dots
\end{aligned}
\end{equation}

If $p\ge n$, the sequence of remainders above is periodic with period $n$, given by the part of the sequence shown above, and the sum of any subsequence of $n$ consecutive remainders is 0. So if $p$ is divisible by $n$, $P_{p,q}(x)$ is divisible by $Q(x)$. If $p$ is not divisible by $n$, then doing first the cancellation as above and then using the symmetry of the sequence of remainders about the remainder equal to 0, the sum of the sequence of remainders equals 0 only if $r+1=t$, that is, only if $q$ is divisible by $n$.

\section{Proof of Proposition~\ref{propo:skewd}}

Consider a $k$-inflated copy of the $L$ $n$-omino. Using the presence of $B_3(k)$ in the Gr\" obner basis, and Theorem~\ref{thm-groebner-tiling}, the existence of a signed tiling of the copy becomes equivalent to deciding when a $k\times nk$ rectangle has a signed tiling by $\mathcal{T}_n$. Theorem~\ref{thm-main} implies that this is always the case.

\section{Proof of Proposition~\ref{propo:skew93}}

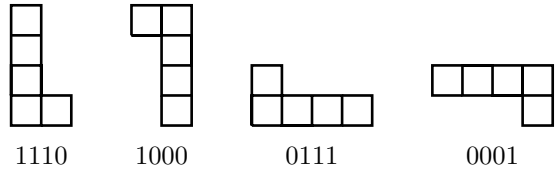
\begin{figure}[h!]
\centering
\begin{tikzpicture}[scale=.4]
\draw [line width = 1] (1, 5)--(2, 5) -- (3, 5) -- (4, 5) -- (4, 6) -- (2, 6) -- (2, 7) -- (1, 7) -- (1,5) -- (1, 6) -- (2, 6) -- (2, 5) -- (3, 5) -- (3, 6);
\draw [line width = 1] (4,5)--(5,5)--(5,6)--(4,6);
\node at (3,4) {0111};
\draw [line width = 1] (8, 7) -- (11, 7) -- (11, 5) -- (10, 5) -- (10, 6) -- (8, 6) -- (8, 7) -- (9, 7) -- (9, 6) -- (10, 6) -- (10, 7) -- (10, 6) -- (11, 6);
\draw [line width = 1] (8,7)--(7,7)--(7,6)--(8,6);
\node at (9,4) {0001};
\draw [line width = 1](-6, 5) -- (-5, 5) -- (-5, 6) -- (-6, 6) -- (-6, 5) -- (-7, 5) -- (-7, 6) -- (-6, 6) -- (-6, 7) -- (-7, 7) -- (-7, 6) -- (-7, 9) --
(-6, 9) -- (-6, 6);
\draw [line width = 1] (-7,8)--(-6,8);
\node  at (-6,4) {1110};
\draw [line width = 1](-3, 8) -- (-3, 9) -- (-1, 9) -- (-1, 6) -- (-2, 6) -- (-2, 8) -- (-3, 8) -- (-2, 8) -- (-2, 9) -- (-2, 8) -- (-1, 8) -- (-1, 7) -- (-2, 7);
\draw [line width = 1] (-2,6)--(-2,5)--(-1,5)--(-1,6);
\node  at (-2,4) {1000};
\end{tikzpicture}
\caption{The four $L$-shaped ribbon pentominoes and their encodings}
\label{f:encodings}
\end{figure}

1) We employ a ribbon tiling invariant introduced by Pak~\cite{Pak}. Each ribbon tile of length $n$ can be encoded uniquely as a binary string of length $n-1,$ denoted $(\epsilon_1, \dots, \epsilon_{n-1}),$ where a $1$ represents a down movement and a $0$ represents a right movement. The encoding of a $1\times n$ bar is $(0,0,\dots,0),$ for a $n\times 1$ bar is $(1,1,\dots,1),$ and for the tiles in $\mathcal{T}_5$ the encodings are shown in Figure~\ref{f:encodings}.

Pak showed that the function $f_1(\epsilon_1, \dots, \epsilon_{n-1}) = \epsilon_1-\epsilon_{n-1}$ is an invariant of the set of ribbon tiles made of $n$-cells, which contains as a subset $\mathcal{T}_n$. In particular, one has that
\begin{equation}\label{eq:inv-tiling-223}
f_1(\epsilon_1, \dots, \epsilon_{n-1}) =\pm 1
\end{equation}
for any tile in $\mathcal{T}_n$. The area of a $k$-inflated copy of the $L$ $n$-omino is an odd multiple of $n$ and can be easily covered by $1\times n$ and $n\times 1$ bars, each one having the invariant equal to zero. If we try to tile by $\mathcal{T}_n$, then the invariant is zero only if we use an even number of tiles. But this is impossible because the area is odd.

\begin{figure}[h!]
\centering
\begin{tikzpicture}[scale=.36]

\draw [line width = 1] (-3, 0)--(5, 0) -- (5, 4) -- (1, 4) -- (1, 5) -- (-3, 5) -- (-3, 0);
\node at (-4.5,3) {$n-r$};
\node at (-1,5.5) {$r$};
\node at (3,4.5) {$r$};
\node at (5.5,2.5) {$r$};

\draw [line width = 1] (9,0)--(29,0)--(29,10)--(19,10)--(19,7)--(9,7)--(9,0);
\node at (7.5,4.5) {$n-r$};
\node at (14,7.5) {$r$};
\node at (23,10.5) {$r$};
\node at (29.5,5.5) {$r$};

\end{tikzpicture}
\caption{Leftover regions.}
\label{fig:leftovers34}
\end{figure}
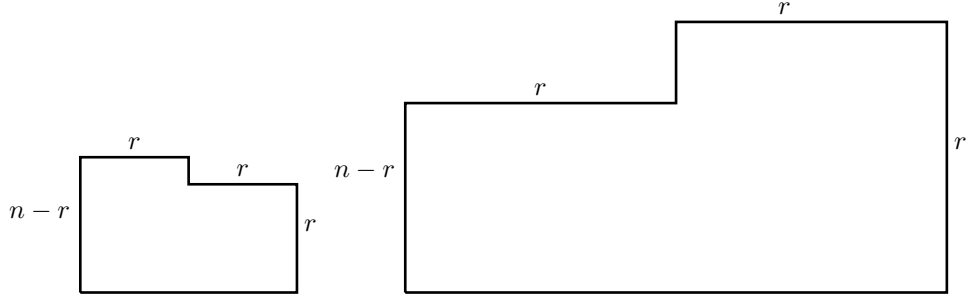

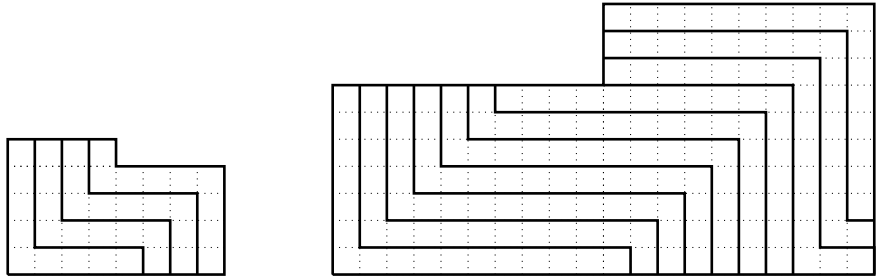
\begin{figure}[h!]
\centering
\begin{tikzpicture}[scale=.36]
\draw[dotted,step=1] (-3, 0) grid (5, 4);
\draw[dotted,step=1] (-3, 4) grid (1, 5);
\draw [line width = 1] (-3, 0)--(5, 0) -- (5, 4) -- (1, 4) -- (1, 5) -- (-3, 5) -- (-3, 0);
\draw [line width = 1] (-2,5)--(-2,1)--(2,1)--(2,0);
\draw [line width = 1] (-1,5) -- (-1, 2) -- (3, 2) -- (3, 0);

\draw [line width = 1] (0, 5) -- (0, 3) -- (4 ,3)--(4, 0);

\draw[dotted,step=1] (9, 0) grid (19, 7);
\draw[dotted,step=1] (19, 0) grid (29, 10);
\draw [line width = 1] (9,0)--(29,0)--(29,10)--(19,10)--(19,7)--(9,7)--(9,0);

\draw [line width = 1] (10,7)--(10,1)--(20,1)--(20,0);
\draw [line width = 1] (11,7)--(11,2)--(21,2)--(21,0);
\draw [line width = 1] (12,7)--(12,3)--(22,3)--(22,0);
\draw [line width = 1] (13,7)--(13,4)--(23,4)--(23,0);
\draw [line width = 1] (14,7)--(14,5)--(24,5)--(24,0);
\draw [line width = 1] (15,7)--(15,6)--(25,6)--(25,0);
\draw [line width = 1] (19,7)--(26,7)--(26,0);
\draw [line width = 1] (19,8)--(27,8)--(27,1)--(29,1)--(29,0);
\draw [line width = 1] (19,9)--(28,9)--(28,2)--(29,2);

\end{tikzpicture}
\caption{Tiling the leftover region by ribbon $n$ tiles, cases $n=5,k=4$, and $n=17, k=20$.}
\label{fig:fine-strips}
\end{figure}

2) Let $k=n\ell+r, 0<r<n.$ After cutting from a $k$-inflated copy a region that can be covered by $1\times n$ and $n\times 1$ bars, and which has the $f_1$ invariant equal to zero, we are left with one of the regions shown in Figure~\ref{fig:leftovers34}. Case a) appears if $2r<n$ and case b) appears if $2r>n$.
Both of these regions can be tiled by $r$ ribbon tiles of area $n$ as in Figure~\ref{fig:fine-strips}. In the first case the sequence of $r$ encodings of the ribbon tiles is:
\begin{equation}
\begin{aligned}
1,1,1,\dots,1,1,1,&0,0,\dots,0,0,0\\
1,1,1,\dots,1,1,0,&0,0,\dots,0,0,1\\
1,1,1,\dots,1,0,0,&0,0,\dots,0,1,1\\
\dots &\dots\dots
\end{aligned}
\end{equation}
where we start with $n-r-1$ ones and $r$ zeros, and then shift the zeroes to the left by 1 at each step, completing the sequence at the end with ones.
As $r\le n-r-1$, the subsequence of zeroes does not reach the left side, so the $f_1$ invariant of the region is equal to $1$.

In the second case, the sequence of $r$ encodings of the ribbon tiles starts as above, but now the subsequence of zeroes reaches the left side. Then we have a jump of $n-r$ units of the sequence of zeroes to the left, the appearance of an extra one at the right, and a completion of the sequence by zeroes to the right. Then the subsequence of ones that appears start shifting to the right till it reaches the right edge. The $f_1$ invariant of the region is equal to $-1$.

So in both cases the $f_1$ invariant is an odd number. Nevertheless, if the $k$-copy is tiled by $\mathcal{T}_n$, one has to use an even number of tiles and the invariant is an even number. Contradiction.

\section{Proof of Theorem~\ref{thm:last63}}

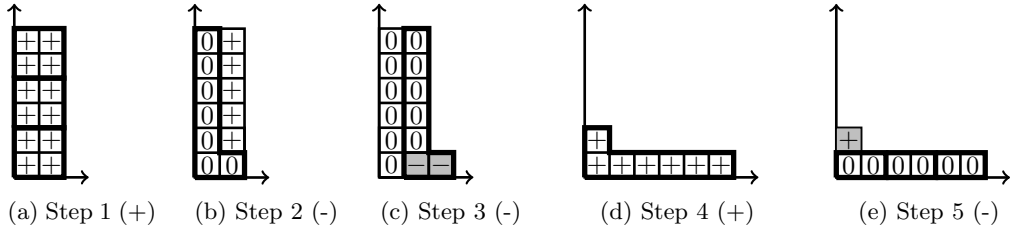
\begin{figure}[h!]
\centering
\begin{subfigure}{.12\textwidth}
\begin{tikzpicture}[scale=.33]

\draw [line width = 1,<->] (0,7)--(0,0)--(3,0);
\draw [line width = 2] (0,0)--(2,0)--(2,2)--(0,2)--(0,0);
\draw [line width = 1] (1,0)--(1,2);
\draw [line width = 1] (0,1)--(2,1);

\node at (.5,.5) {$+$};
\node at (.5,1.5) {$+$};
\node at (1.5,.5) {$+$};
\node at (1.5,1.5) {$+$};

\draw [line width = 2] (0,2)--(2,2)--(2,4)--(0,4)--(0,2);
\draw [line width = 1] (1,2)--(1,4);
\draw [line width = 1] (0,3)--(2,3);

\node at (.5,2.5) {$+$};
\node at (.5,3.5) {$+$};
\node at (1.5,2.5) {$+$};
\node at (1.5,3.5) {$+$};

\draw [line width = 2] (0,4)--(2,4)--(2,6)--(0,6)--(0,4);
\draw [line width = 1] (1,4)--(1,6);
\draw [line width = 1] (0,5)--(2,5);

\node at (.5,4.5) {$+$};
\node at (.5,5.5) {$+$};
\node at (1.5,4.5) {$+$};
\node at (1.5,5.5) {$+$};
\end{tikzpicture}
\caption{Step 1 (+)}
\end{subfigure}
~~
\begin{subfigure}{.12\textwidth}
\begin{tikzpicture}[scale=.33]

\draw [line width = 1,<->] (0,7)--(0,0)--(3,0);
\draw [line width = 1] (0,0)--(2,0)--(2,2)--(0,2)--(0,0);
\draw [line width = 2] (0,0)--(2,0)--(2,1)--(1,1)--(1,6)--(0,6)--(0,0);
\draw [line width = 1] (1,0)--(1,2);
\draw [line width = 1] (0,1)--(2,1);

\draw [line width = 1] (0,5)--(1,5);
\draw [line width = 1] (0,4)--(1,4);
\draw [line width = 1] (0,3)--(1,3);
\draw [line width=1] (2,1)--(2,6)--(1,6);

\draw [line width = 1] (1,5)--(2,5);
\draw [line width = 1] (1,4)--(2,4);
\draw [line width = 1] (1,3)--(2,3);

\node at (.5,2.5) {$0$};
\node at (.5,3.5) {$0$};
\node at (.5,4.5) {$0$};
\node at (.5,5.5) {$0$};

\node at (.5,.5) {$0$};
\node at (.5,1.5) {$0$};
\node at (1.5,.5) {$0$};
\node at (1.5,1.5) {$+$};
\node at (1.5,2.5) {$+$};
\node at (1.5,3.5) {$+$};
\node at (1.5,4.5) {$+$};
\node at (1.5,5.5) {$+$};
\end{tikzpicture}
\caption{Step 2 (-)}
\end{subfigure}
~~
\begin{subfigure}{.12\textwidth}
\begin{tikzpicture}[scale=.33]

\draw [line width = 1,<->] (0,7)--(0,0)--(3.8,0);
\draw [fill=lightgray] (1,0)--(3,0)--(3,1)--(1,1)--(1,0);

\draw [line width = 1] (0,0)--(2,0)--(2,2)--(0,2)--(0,0);
\draw [line width = 1] (0,0)--(2,0)--(2,1)--(1,1)--(1,6)--(0,6)--(0,0);
\draw [line width = 2] (1,0)--(3,0)--(3,1)--(2,1)--(2,6)--(1,6)--(1,0);
\draw [line width = 1] (1,0)--(1,2);
\draw [line width = 1] (0,1)--(2,1);

\draw [line width = 1] (0,5)--(1,5);
\draw [line width = 1] (0,4)--(1,4);
\draw [line width = 1] (0,3)--(1,3);
\draw [line width=1] (2,1)--(2,6)--(1,6);

\draw [line width = 1] (1,5)--(2,5);
\draw [line width = 1] (1,4)--(2,4);
\draw [line width = 1] (1,3)--(2,3);

\node at (.5,2.5) {$0$};
\node at (.5,3.5) {$0$};
\node at (.5,4.5) {$0$};
\node at (.5,5.5) {$0$};

\node at (.5,.5) {$0$};
\node at (.5,1.5) {$0$};
\node at (1.5,.5) {$-$};
\node at (1.5,1.5) {$0$};
\node at (1.5,2.5) {$0$};
\node at (1.5,3.5) {$0$};
\node at (1.5,4.5) {$0$};
\node at (1.5,5.5) {$0$};
\node at (2.5,.5) {$-$};
\end{tikzpicture}
\caption{Step 3 (-)}
\end{subfigure}
~~~~
\begin{subfigure}{.16\textwidth}
\begin{tikzpicture}[scale=.33]

\draw [line width = 1,<->] (0,7)--(0,0)--(7,0);
\draw [line width = 2] (0,0)--(0,2)--(1,2)--(1,1)--(6,1)--(6,0)--(0,0);

\draw [line width = 1] (1,0)--(1,1);
\draw [line width = 1] (0,1)--(1,1);
\draw [line width = 1] (2,0)--(2,1);
\draw [line width = 1] (3,0)--(3,1);
\draw [line width = 1] (4,0)--(4,1);
\draw [line width = 1] (5,0)--(5,1);

\node at (.5,1.5) {$+$};
\node at (.5,.5) {$+$};
\node at (1.5,.5) {$+$};
\node at (2.5,.5) {$+$};
\node at (3.5,.5) {$+$};
\node at (4.5,.5) {$+$};
\node at (5.5,.5) {$+$};

\end{tikzpicture}
\caption{Step 4 (+)}
\end{subfigure}
~~~~
\begin{subfigure}{.16\textwidth}
\begin{tikzpicture}[scale=.33]

\draw [line width = 1,<->] (0,7)--(0,0)--(7,0);
\draw [line width = 1] (0,0)--(0,2)--(1,2)--(1,1)--(6,1)--(6,0)--(0,0);

\draw [fill=lightgray] (0,1)--(1,1)--(1,2)--(0,2)--(0,1);

\draw [line width = 2] (0,0)--(2,0)--(2,1)--(0,1)--(0,0);
\draw [line width = 2] (2,0)--(4,0)--(4,1)--(2,1);
\draw [line width = 2] (4,0)--(6,0)--(6,1)--(4,1);

\draw [line width = 1] (1,0)--(1,1);
\draw [line width = 1] (0,1)--(1,1);
\draw [line width = 1] (2,0)--(2,1);
\draw [line width = 1] (3,0)--(3,1);
\draw [line width = 1] (4,0)--(4,1);
\draw [line width = 1] (5,0)--(5,1);

\node at (.5,1.5) {$+$};
\node at (.5,.5) {$0$};
\node at (1.5,.5) {$0$};
\node at (2.5,.5) {$0$};
\node at (3.5,.5) {$0$};
\node at (4.5,.5) {$0$};
\node at (5.5,.5) {$0$};

\end{tikzpicture}
\caption{Step 5 (-)}
\end{subfigure}
~~~~
\caption{Generating a single cell by $\tilde{\mathcal{T}}_n$.}
\label{fig:null392}
\end{figure}

It is enough to generate the tile consisting of a single cell. We show the proof for $n=7$ in Figure~\ref{fig:null392}. The proof can be easily generalized to any $n\ge 5,$ odd. First we construct a domino with both cells having the same sign (as in Figure~\ref{fig:null392} c)), and then we use it to reduce  the $L$ $n$-omino until a single cell is left.

\section{The method of Barnes}

In this section we give a proof of Theorem~\ref{thm-main-barnes} following a method developed by Barnes. The reader of this section should be familiar with \cite{Barnes1, Barnes2}. We apply the method to the infinite collection of tiling sets $\mathcal{T}_n, n\ge 5$ odd.

Let $n\ge 5$ odd fixed. Consider the polynomials  \eqref{eq:generators-k} associated to the tiles in $\mathcal{T}_n$ and denote by $I$ the ideal generated $G_1(k)$, $G_2(k)$, $G_3(k)$, $G_4(k)$. We show that the algebraic variety $V\subset \mathbb{C}^2$ defined by  \eqref{eq:generators-k} is zero dimensional and consists only of the pairs of points \begin{equation}
\left (\epsilon, \frac{1-\epsilon^n}{1-\epsilon}\right ),
\end{equation}
where $\epsilon$ is an $n$-th root of identity different from 1.

Separate $x$ from $G_1(k)=0$ and replace in $G_2(k)=0$ to have:
\begin{equation}\label{eq:points-v}
y^{2k-1}-\frac{y^{2k}-1}{y-1}\cdot\frac{y^{2k}-1}{y-1}=0.
\end{equation}

Eliminating the denominators gives:
\begin{equation}
y^{2k-1}-(y^{2k-1}+y^{2k-2}+\dots +y^2+y+1)^2=0,
\end{equation}
which can be factored as:
\begin{equation}
(y^{2k}+y^{2k-1}+y^{2k-2}+\dots +y^2+y+1)(y^{2k-2}+y^{2k-3}+\dots +y^2+y+1)=0.
\end{equation}

It is clear that all roots of the polynomial above, and of the corresponding polynomial in the variable $x$, are roots of unity of order $2k+1$ and $2k-1$. Using the system of equations that defines $V$, the roots of order $2k-1$ can be eliminated. Moreover, the only solutions of the system are given by \eqref{eq:points-v}.

We show now that $I$ is a radical ideal. For this we use an algorithm of Seidenberg which can be applied to find the radical ideal of a zero dimensional algebraic variety over an algebraically closed field. See Lemma 92 in~\cite{Seidenberg}. Compare also with \cite[Theorem 7.1]{Barnes1}. As $V\subset \mathbb{C}^2$ is zero dimensional, one can find $f_1(x)$ and $f_2(y)$ that belong to the radical ideal. We consider the square free polynomials:
\begin{equation*}
\begin{gathered}
f_1(x)=x^{2k}+x^{2k-1}+x^{2k-2}+\dots +x^2+x+1,\\
f_2(y)=y^{2k}+y^{2k-1}+y^{2k-2}+\dots +y^2+y+1.
\end{gathered}
\end{equation*}

If $f_1(x), f_2(y)$ are square free, then the ideal generated by $G_1(k)$, $G_2(k)$, $G_3(k)$, $G_4(k)$ and $f_1(x), f_2(y)$ is radical. So, in order to show that $I$ itself is radical, it is enough to show that $f_1(x), f_2(y)$ belong to $I$.
It is easier to generate $f_i$ using the Gr\"obner basis, so we will use this approach.

\begin{propo} The polynomials $f_1(x), f_2(y)$ belong to the ideal $I$.
\end{propo}

\begin{pf} It is enough to generate $f_1(x)$. One has:
\begin{equation*}
f_1(x)=xG_3(x)-B_3(k).
\end{equation*}
\end{pf}

We can apply now the main result in \cite[Lemma 3.8]{Barnes1}: a region $R$ is signed tiled by $\mathcal{T}_n$ if and only if the polynomial $f_R(x,y)$ associated to $R$ evaluates to zero in any point of the variety $V$. If $R$ is a rectangle of dimensions $p\times q$ in the square lattice, then
\begin{equation}
f_R(x,y)=\frac{x^q-1}{x-1}\cdot \frac{y^p-1}{y-1},
\end{equation}
which clearly evaluates to zero in all points of $V$ if and only if one of $p, q$ is divisible by $n$.

The fact that Theorem~\ref{thm-main-barnes} implies Theorem~\ref{thm-main} follows the idea of \cite[Theorem 4.2]{Barnes1}. Indeed, a set of generators for the regions that are signed tiled with $\mathcal{Q}$ by $\mathcal{T}_n$ is given by the polynomials $f_1(x), f_2(y)$ above. Both of them can be generated by the Gr\"obner basis using only integer coefficients.

\section*{Acknowledgement}

The author was partially supported by Simons Foundation Grant 208729.

\end{document}